\documentclass[a4paper,12pt,twoside]{article}
\usepackage{amssymb}
\usepackage[final]{lpretex}
\usepackage{title}
\usepackage{a4}
\usepackage{amscd}
\input xy
\xyoption{all}
\usepackage{graphicx}
\EndOfDump

\def\DHpartformat#1{(#1) }
\def\DHrefpart#1{(\DHRefpart{#1})}
\LBsetstyleof{partno}{arabic}

\def\i {\item}

\let\define\def

\def\A {{\mathbb A}} \def\B {{\mathbb B}} \def\C {{\mathbb C}}
\def\D {\Delta}  \def\F {{\mathbb F}}

\def\GG {{\mathbb G}}   
  
  \def\P {{\mathbb P}} 
\def\Q {{\mathbb Q}} \def\R {{\mathbb R}}
 
  \def\X {{\mathbb X}}
\def\Z {{\mathbb Z}} 

\define \n {\mathbb N}
\define \z {\mathbb Z}
\define \q {\mathbb Q}
\define \PP {\mathbb P}

\def\sA {{\Cal A}} \def\sB {{\Cal B}} \def\sC {{\Cal C}}
\def\sD {{\Cal D}} \def\sE {{\Cal E}} \def\sF {{\Cal F}}
\def\sG {{\Cal G}} \def\sH {{\Cal H}} \def\sI {{\Cal I}}
  \def\sL {{\Cal L}}
\def\sM {{\Cal M}} \def\sN {{\Cal N}} \def\sO {{\Cal O}}
\def\sQ {{\Cal Q}} \def\sR {{\Cal R}}
 \def\sT {{\Cal T}} \def\sU {{\Cal U}}
\def\sV {{\Cal V}} \def\sW {{\Cal W}} \def\sX {{\Cal X}}
\def\sY {{\Cal Y}}

\def\sZ {{\Cal Z}}
\define \cN {\Cal N}
\define \cf {\Cal F}
\define \cg {\Cal G}
\define \cE {\Cal E}
\define \ce {\Cal E}
\define \cc {\Cal C}
\define \cV {\Cal V}
\define \cA {\Cal A}
\define \cK {\Cal K}
\define \cO {\Cal O}
\define \cF {\Cal F}
\define \cn {\Cal N}
\define \cI {\Cal I}
\define \sP {\Cal P}
\def\a {\alpha} \def\b {\beta} \def\g {\gamma} \def \d {\delta} 
\def\e {\epsilon}   
\def\s {\sigma}

\define \x {\xi}
\define \y {\eta}
\define \G {\Gamma}
\define \r {\rho}
\define \w {\omega}

\def \tZ {\widetilde Z}
\def\tX {\widetilde X}

\def \trho {\tilde {\rho}}

\def \tpi {\tilde{\pi}}

\define \tH {\widetilde H}
\define \tG {\widetilde{G}}
\define \tW {\widetilde W}
\define \tF {\widetilde F}
\define \tm {\tilde m}
\define \St {\widetilde S}
\define \Xt {\widetilde X}
\define \tS {\widetilde S}
\define \tpsi {\tilde \psi}
\define \tL {\widetilde L}
\define \tE {\widetilde E}
\define \tl {\tilde l}
\define \tA {\widetilde A}
\define \tom {\tilde\omega}
\define \tT {\widetilde T}
\define \tB {\widetilde B}
\define \tf {\tilde f}
\define \tsA {\widetilde{\sA}}

\define \tM {\widetilde M}
\define \tsM {\widetilde{\sM}}
\define \tpsi {\widetilde{\psi}}
\define \trho {\widetilde{\rho}}
\define \tR {\widetilde R}
\define \tp {\widetilde p}
\define \tq {\widetilde q}
\define \tc {\tilde c}
\define \tsF {\widetilde {\sF}}
\define \tx {\tilde x}
\define \tg {\tilde g}
\define \tw {\tilde w}
\define \ta {\widetilde\alpha}
\define \tsZ {\widetilde{\sZ}}
\define\ts{\tilde\sigma}

\def\pd {\partial}

\def \Dx1 {\frac{\pd}{{\pd} x_1}}
\def \Dy1 {\frac{\pd}{{\pd} y_1}}
\def \Dz1 {\frac{\pd}{{\pd} z_1}}
\def \Dx2 {\frac{\pd}{{\pd} x_2}}
\def \Dy2 {\frac{\pd}{{\pd} y_2}}
\def \Dz2 {\frac{\pd}{{\pd} z_2}}

\def\q {\quad} 

\def\Mapdiagr#1{\nearrow\rlap{$\lower 5pt\vbox{{\hbox{$\mkern
-15mu\scriptstyle#1$}}}$}} 

\def\Mapdiagl#1{\llap{$\lower 5pt\vbox{{\hbox{$\scriptstyle#1\mkern
-15mu$}}}$}\searrow} 

\def\Mapswr#1{\swarrow\rlap{$\lower 5pt\vbox{{\hbox{$\mkern
-15mu\scriptstyle#1$}}}$}}              

\def\Mapnwl#1{\nwarrow\rlap{$\lower 5pt\vbox{{\hbox{$\mkern
-15mu\scriptstyle#1$}}}$}}

\def\i.e#1#2#3{\mathrel{\smash{\mathop{#2}\limits^{#1}_{#3}}}}

\def \inj {\hookrightarrow}
\def \onto {\twoheadrightarrow}  

\define \Rhook {\hookrightarrow}

\def \half {\raise1pt\hbox{$\scriptstyle
        \frac{1}{2}\displaystyle$}}

\def \x{{\sl X}\llap{$\mkern -2mu {\scriptstyle -}$}}

\def \Hom {\operatorname{Hom}}

\def \Proj {\operatorname{Proj}}
\def \Symm {\operatorname{Sym}}

\def \Bl {\operatorname{Bl}}
\def \Pic {\operatorname{Pic}}

\define \Kod {\operatorname{Kod}}
\define \dimension {\operatorname{dim}}
\define \codim {\operatorname{codim}}
\define \contr {\operatorname{contr}}
\define \rk {\operatorname{rank}}
\define \Im {\operatorname {Im}}
\define \Mor {\operatorname{Mor}}
\define \Cl {\operatorname{Cl}}
\define \Hilb {\operatorname{Hilb}}
\define \degree {\operatorname{deg}}
\define \mult {\operatorname{mult}}
\define \Aut {\operatorname{Aut}}
\define \NS {\operatorname{NS}}
\define \Gal {\operatorname{Gal}}
\define \ch {\operatorname{char}}
\define \Jac {\operatorname{Jac}}
\define \Km {\operatorname{Km}}
\define \Sec {\operatorname{Sec}}
\define \Stab {\operatorname{Stab}}
\define \Br {\operatorname{Br}}
\define \Inv {\operatorname {Inv}}
\define \tr {\operatorname{tr}}
\define \Frob {\operatorname{Frob}}
\define \Symn {\operatorname{Sym}^n}
\define \Ev {\sE^\vee}
\define \ordp {\operatorname{ord}_p}
\define \Supp {\operatorname{Supp}}
\define \Ann {\operatorname{Ann}}
\define \disc {\operatorname{disc}}
\define \lie {\operatorname{lie}}
\define \embdim {\operatorname{embdim}}

\def\ad{\operatorname{ad}}

\def\Ob{\operatorname{Ob}}

\define\bark{\overline{k}}
\define\Lie{\operatorname{Lie}}
\define\tsG{\widetilde{\sG}}
\def\uu{\underline{\mathfrak u}}

\def\gg{\underline{\mathfrak g}}

\def\ll{\underline{\mathfrak l}}
\def\nn{\underline{\mathfrak n}}
\define\tsH{\widetilde{\sH}}
\define\hY{\widehat{Y}}

\define\szd{\sZ^{\dagger}}

\define\tszd{{\widetilde{\sZ}}^{\dagger}}
\define\bX{\overline{X}}
\define\tsW{\widetilde{\sW}}
\define\tsE{\widetilde{\sE}}

\define\hsY{\widehat{\sY}}

\define\szdss{\sZ^{\dagger ss}}

\define\barpsi{\overline{\psi}}
\define\bsF{\overline{\sF}}
\define\bsH{\overline{\sH}}
\define\bF{\overline{F}}
\define\bH{\overline{H}}

\define\tphi{\widetilde{\phi}}

\def\hod#1#2#3#4{\ensuremath{ if#30 H^{#2}({#1},{\cal O}_{#1}) \else 
 H^{#2}(#1,\Omega^{#3}\if\relax{#4}\relax_{#1}\else _{#1/#4}\fi)\fi}}

\begin{document}
\title
[del Pezzo surfaces as Springer fibres]
{Del Pezzo surfaces as Springer fibres for exceptional groups}
\author{I. Grojnowski}
\address{D.P.M.M.S.\\
Cambridge University\\
Cambridge CB2 1SB\\
U.K.}
\email{groj@dpmms.cam.ac.uk}
\author{N. I. Shepherd-Barron}
\address{King's College London\\
Strand\\
London WC2R 2LS\\
U.K.}
\email{Nicholas.Shepherd-Barron@kcl.ac.uk}
\maketitle

\begin{section}{Introduction}
Since the discovery by Cayley and Salmon of the $27$ lines on a cubic surface
and their configuration 
and that by Killing and Cartan of the exceptional 
simple Lie algebra
of type $E_6$ it has been clear that del Pezzo surfaces and
exceptional simple algebraic groups have attached to them
the same combinatorial objects.
That is, there are well known constructions
$${\mathrm{(del\ Pezzo\ surfaces)}}\rightarrow
{\mathrm{(exceptional\ root\ data)}}\leftrightarrow
{\mathrm{(exceptional\ simple\ groups)}}.$$

One of the results of this paper is a direct geometrical
construction, for $5\le l\le 8$, of the del Pezzo surfaces $S$ of 
degree $d=9-l$
from the split simply connected exceptional simple\footnote
{We write $E_5=D_5$ and regard $D_5$ as exceptional.
We also refer to groups as ``simple'' rather than
``almost simple''.}
group $G$ of type $E_l$, in the presence of the universal
elliptic curve $\sE$. Direct, that is, in the sense
that it does not involve root data. 

The del Pezzo surfaces arise out of the group as follows. Let $P$ denote the weight
lattice of $G$ and $B$ a Borel subgroup.
Suppose that $\tsG^{ss}_{\sE}$ is the stack that classifies
semi-stable $G$-bundles over elliptic curves together with
a reduction of degree $0$ of the structure group to $B$.
(The existence of such a reduction is a characterization
of semi-stability.)
Then there is a partial relative compactification 
of $\tsG^{ss}_{\sE}$
in which the total space $\sD_1^-$ of 
a certain family $\sD_1^-\to \sY=\Hom(P,\sE)$
of these del Pezzo surfaces appears
as a boundary divisor.
We shall describe this in more detail later in this introduction.
(The del Pezzo surfaces in the family $\sD_1^-\to \sY$
are in fact \emph{marked}; that is, 
a copy of $\sE$ is embedded in them as an
anticanonical curve and their N{\'e}ron--Severi groups
are rigidified. The family is nearly, but not
quite, universal; the fine moduli space
is a non-separated union of copies of $\sY$,
and $\sY$ is then its maximal separated quotient.)

However, we have no suggestion for
a similarly direct construction of the groups from the 
del Pezzo surfaces
that would bypass consideration of the root data.

The title of this paper refers to
our presentation of these results as an extension of the 
construction and theorems
due to Brieskorn, Grothendieck, Slodowy and Springer
(BGSS) and described in detail in [Sl] 
that, in good characteristics, 
reveal the deformations and
simultaneous resolution of du Val singularities
(= rational double points = simple singularities = 
Kleinian singularities)
of type $A,D$ or $E$ inside
the corresponding split simply connected simple algebraic group. 
In turn, this provides a characterization of du Val
singularities. Of course, there are many others.
One, that we shall use, is that a surface
singularity is du Val if and only if
its monodromy group, acting
on the vanishing cohomology, is the relevant Weyl group.

Our results lead, as a by-product, to an identification of the unipotent 
singularity of $E_8$ in bad characteristics 
(that is, $2,3$ and $5$).
This is done in section \ref{weird}.
The proof depends upon specializing
$\sE$ to a uniformizable
elliptic curve $E$ over a complete algebraically closed field
of characteristic $p$. 
The picture is that the curve $E$
acts as a vector to carry some of the geometry
of the group into a del Pezzo surface $S$
with a singularity isomorphic to the unipotent singularity,
and we can then exploit the fact that
the defining equations of such a surface
are so simple that classification
of the unipotent singularities becomes very easy and quick.
Moreover, this leads to an embedding
of the group-theoretical description of 
the deformations and resolutions of du Val singularities
into their description in terms of del Pezzo surfaces,
as expounded in [SSS].

This is, perhaps, some slight evidence
supporting a view that our results 
tell us something about the group, and that further information
about the group is to be found in
higher codimension in suitable relative compactifications 
of $\tsG^{ss}_{\sE}$.

In the rest of this introduction
we give a more detailed survey of our results.
 
From the viewpoint of the geometry of algebraic
surfaces and their degenerations
the extension referred to above consists of showing that
{\emph{simultaneous log resolutions}},
which are defined in the next paragraph, of
simply elliptic singularities of degrees $d=1,2,3$ or $4$,
which were first observed in the context of type II
degenerations of K3 surfaces, 
although inspired by work of 
Brieskorn [B] and Looijenga [L],
also have realizations in terms of simple algebraic groups.
Here, only exceptional groups of type $E_{9-d}$ can, and do, occur.

We define a \emph{del Pezzo surface}
to be a surface $S$ whose anti-canonical divisor class is ample 
and that has at worst
du Val singularities; a \emph{weak del Pezzo surface}
is smooth and its anti-canonical class is ample modulo finitely many
$(-2)$-curves. That is, a weak del Pezzo surface is the minimal
resolution of a del Pezzo surface
and a del Pezzo surface is the anticanonical model of a weak
del Pezzo surface.

\begin{definition}\label{log} Given a flat family $X\to S$ of 
normal Gorenstein
surfaces, a {\emph{simultaneous log resolution}} of $X\to S$ is a 
commutative diagram
$$\xymatrix{
{\tX}\ar[r]\ar[d] & {X}\ar[d]\\
{\tS}\ar[r] & {S}
}$$
where $\tS\to S$ is proper, dominant
and generically finite, $\tX\to X\times_S\tS$ is proper and birational,
$\tX$ is smooth, $\tX\to\tS$ is flat
and semi-stable, in the sense that each geometric fibre
is a reduced union of smooth surfaces with normal crossings, 
and the relative canonical class $K_{\tX/\tS}$ is 
equivalent to the pullback of $K_{X/S}$.
\end{definition}
It is well known in the context of
the birational geometry of degenerating families of K3 surfaces
that a versal deformation $X\to S$
of a simply elliptic singularity 
of degree at most $4$
possesses a simultaneous log resolution $\tX\to\tS$ 
with the following properties.

\noindent (1) $\tX\to\tS$ is \emph{of type II}, in the sense that,
for each fibre $X_s$ with an elliptic singularity,
its inverse image $\tX_{\tilde s}$ in $\tX$ has two components,
say $\tX_{\tilde s}=D_1+D_2$.
Moreover $D_1$ 
is the minimal resolution of $X_s$, $D_2$ is a 
weak del Pezzo surface, $D_2$ is contracted to a point 
under $\tX\to X$ and $D_1\cap D_2$ is the elliptic curve that is the
exceptional locus of $D_1\to X_s$.

\noindent (2) For each fibre $X_s$ that has at most du Val singularities,
its inverse image $\tX_{\tilde s}$ in $\tX$ is the minimal
resolution of $X_s$.

\noindent (3) the base change $\tS\to S$ is the composite of
a ramified Galois covering of the base $S$ 
whose Galois group is the corresponding finite Weyl group and 
a blow-up along the simply elliptic locus in $S$.

On the other hand, Helmke and Slodowy [HS1] and [HS3]
have given a brief indication
of a proof that miniversal deformations of simply elliptic singularities of degree
$d$ can be realized inside $\sG$ when $G=E_{9-d}$.

Recall from [HS2], Theorems $5.6$ and $5.12$, 
that a \emph{regular} 
bundle on an elliptic curve $E$ over a field
is one whose automorphism group has minimal dimension,
namely $l$, while for a \emph{subregular} bundle the automorphism
group has the next smallest dimension, $l+2$. 
Up to translation by points of $E$ there is a unique unstable regular bundle $\eta$
and a unique unstable subregular bundle $\xi$; in a miniversal
deformation of $\xi$ the unstable locus is,
up to this translation, a surface $S$ and all
points of $S-\{\xi\}$ correspond to bundles isomorphic to $\eta$.
Note that in order to construct these bundles it is necessary
(and sufficient) for the group to
be split and the curve $E$ to have a point; it cannot be replaced by
an arbitrary curve of genus $1$.
Then both $\xi$ and $\eta$ are constructed as bundles on $\sE$.

Let's list the objects that appear in our main theorem. These are:

\noindent (1) the stack $\sG=\sG_{\sE}$ of principal
$G$-bundles over the universal elliptic curve $\sE$
(defined over the moduli stack $\sM_{Ell}$ of elliptic curves)
and the open substack $\sG^{ss}$ of semi-stable bundles;

\noindent (2) the stack $\tsG^{ss}$ 
that classifies semi-stable $G$-bundles
together with a reduction of the structure group to a Borel subgroup
such that the associated torus bundle is of degree $0$
(in the course of reaching our main result
we shall show that $\tsG^{ss}$ is proper,
representable and generically finite, of degree equal to 
the order of the Weyl group $W$ of $G$, over $\sG^{ss}$, 
and that, when we restrict attention
to a uniformizable elliptic curve, 
this structure reproduces the BGSS picture); 

\noindent (3) stacks $\tsG$ 
and $\tsG^+$ 
(to be described in Section \ref{G+}),
a projective birational morphism
$\tsG\to\tsG^+$ and a morphism which is
proper but not representable 
(although its fibres have finite automorphism
groups)
$\tsG^+\to\sG$
such that $\tsG\to\sG$ and $\tsG^+\to\sG$
are relative compactifications of $\tsG^{ss}\to\sG^{ss}$; 

\noindent (4) the abelian variety $\sY=\Hom(P,\sE)$ 
over $\sM_{Ell}$ and a certain ample $W$-linearized line
bundle $\sL\to\sY$, which gives a cone $\hsY$ over $\sM_{Ell}$
by contracting the zero section
$0_\sL$ of $\sL$ such that the geometric quotient $[\hsY/W]$ of
$\hsY$ by $W$, relative to
$\sM_{Ell}$, is, by [L], an affine space bundle of rank $l+1$ over $\sM_{Ell}$;

\noindent (5) the unstable subregular bundle $\xi$ mentioned above,
that is defined over $\sM_{Ell}$ and
is unique modulo translation by $\sE$;

\noindent (6) the base $\sZ$ of a deformation
that is minimally versal, modulo the translation by $\sE$ mentioned above,
of $\xi$;

\noindent (7) the products
$\tsZ=\sZ\times_{\sG}\tsG$,
and $\tsZ^+=\sZ\times_{\sG}\tsG^+$
and a birational map (not a morphism) 
$\tsZ^+ -\to\tsZ^-$ which is
constructed as a succession of flops
relative to $\sZ$ 
such that the morphisms
$\tsZ\to\sZ$ and $\tsZ^{\pm}\to\sZ$ are representable.

Moreover, let ${}^0\sG$ be the open substack of
$\sG$ that is the complement of the locus of those
unstable bundles $\Xi$ such that $\dim\Aut(\Xi)>l+4$.
That is, the locus of unstable bundles that are ``worse than
subregular'' has been deleted. In particular,
$\sG^{ss}$ is contained in ${}^0\sG$.
Set ${}^0\tsG=\tsG\times_{\sG}{}^0\sG$
and ${}^0\tsG^+=\tsG^+\times_{\sG}{}^0\sG$, so that
$\sZ$ is, up to turning off $\Pic^1(\sE)$, a chart
in ${}^0\sG$ and
$\tsZ$ is a chart in ${}^0\tsG$.
Then the birational morphism 
$\tsZ\to\tsZ^+$ and the birational maps
$\tsZ^+ -\to\tsZ^-$ and $\tsZ-\to\tsZ^-$
are all isomorphisms over $\sG^{ss}$;
that is, all birational modifications
are made inside the boundary of ${}^0\tsG$
and the commutative square 
$$\xymatrix{
{\tsZ}\ar[r]\ar@{.>}[d]&{\tsZ^+}\ar@{.>}[dl]\ar[d]\\
{\tsZ^-}\ar[r]&{\sZ}
}$$
of rational maps and morphisms is a local description
of a commutative square
$$\xymatrix{
{{}^0\tsG}\ar[r]\ar@{.>}[d]&{{}^0\tsG^+}\ar@{.>}[dl]\ar[d]\\
{{}^0\tsG^-}\ar[r]&{{}^0\sG}
}$$
of rational maps and morphisms between stacks.

What is missing from this picture is a
suitable stack $\tsG^-$ of which ${}^0\tsG^-$
is a natural open substack.

The result of Helmke and Slodowy concerning elliptic
singularities that we mentioned above is that there is a flat
morphism $\sZ\to [\hsY/W]$ that is a minimally versal deformation
of a simply elliptic singularity of degree $d$.
We give a detailed proof of this as part of our main result,
which is that there is also an analogue
of the BGSS construction for the situation
involving $\sE$, as in the following theorem,
which summarizes the results of the paper.
The main result is given more concisely in Theorem \ref{main}.

\begin{theorem}
Suppose that $G=E_{9-d}$. 

\part[i] Let $\tsZ^\sharp$ denote any one of $\tsZ,\tsZ^\pm$.
Then there is a commutative diagram
$$\xymatrix{
{\tsZ^\sharp}\ar[r]\ar[d]^{\pi^\sharp} & {\sZ}\ar[d]\\
{\sL}\ar[r] & {[\hsY/W]}
}$$
where $\pi^\sharp:\tsZ^\sharp\to\sL$ is semi-stable
and smooth over $\sL-0_{\sL}$, the complement of the zero
section $0_{\sL}$.

\part[ii] Suppose that $0$ is a point of $0_{\sL}$ lying over
the origin $0_{\sY}$ of $\sY$. Then
$\tsZ_0=D_0+D_1+Q$, $\tsZ^\pm_0=D_0^\pm+ D_1^\pm$,
each $D_0^\sharp\cap D_1^\sharp$ is a copy of $\sE$,
$\sZ_0$ is a cone over $\sE$ of degree $d$ and 
$\tsZ^-$ is its minimal resolution.
$Q$ is a copy of $\P^1\times\P^1$
and the morphism $\tsZ\to\tsZ^+$ contracts $Q$ to a curve
via a projection $Q\to\P^1$.

\part[iii] When $\tsZ^\sharp=\tsZ^-$ then the diagram
is a type II simultaneous log resolution of $\sZ\to [\hsY/W]$.

\part[iv] The rational map $\tsZ^+ -\to\tsZ^-$
is a sequence of flops and the centre of each flop is the total space of
a smooth family of rational curves over $0_{\sL}$.

\part[v] 
$\sZ\to [\hsY/W]$ is a miniversal deformation of a simply
elliptic singularity of degree $d$ over $\Sp\Q$,
and over $\Sp\Z[1/d]$ if $d\le 3$.

\part[vi]
The morphisms $\tsZ^\sharp\to\sZ\times_{[\hsY/W]}\sL$ are birational.

\part[vii]\label{fibres} The exceptional locus $\sD_1^-$ 
in $\sZ^-$ is the total space
of a family $\sD_1^-\to\sY=0_{\sL}$ of weak del Pezzo surfaces
of degree $d$.
In particular, the fibre over a point $(y,\xi)$ of the morphism
$\sZ^-\to\sZ\times_{[\hsY/W]}\sL$, where $y$ is a point in the
zero section $0_\sL$ of $\sL$, is a weak del Pezzo surface of degree
$d$.

\part[viii] Locally on $\sM_{Ell}$ there is an isomorphism
$\tsZ^-\to\omega_{\sD_1^-/\sY}$.
\noproof
\end{theorem}
The fibres of \ref{fibres} are the Springer fibres of the title.

Figure \ref{fig:fibres Z etc.} 
below shows the structure of the fibres
$\tsZ_0,\tsZ^\pm_0$ and $\sZ_0$.

\begin{figure*}[htbp]
\centerline
{\includegraphics[height=0.5\textheight, width=0.5\textwidth, angle=270]{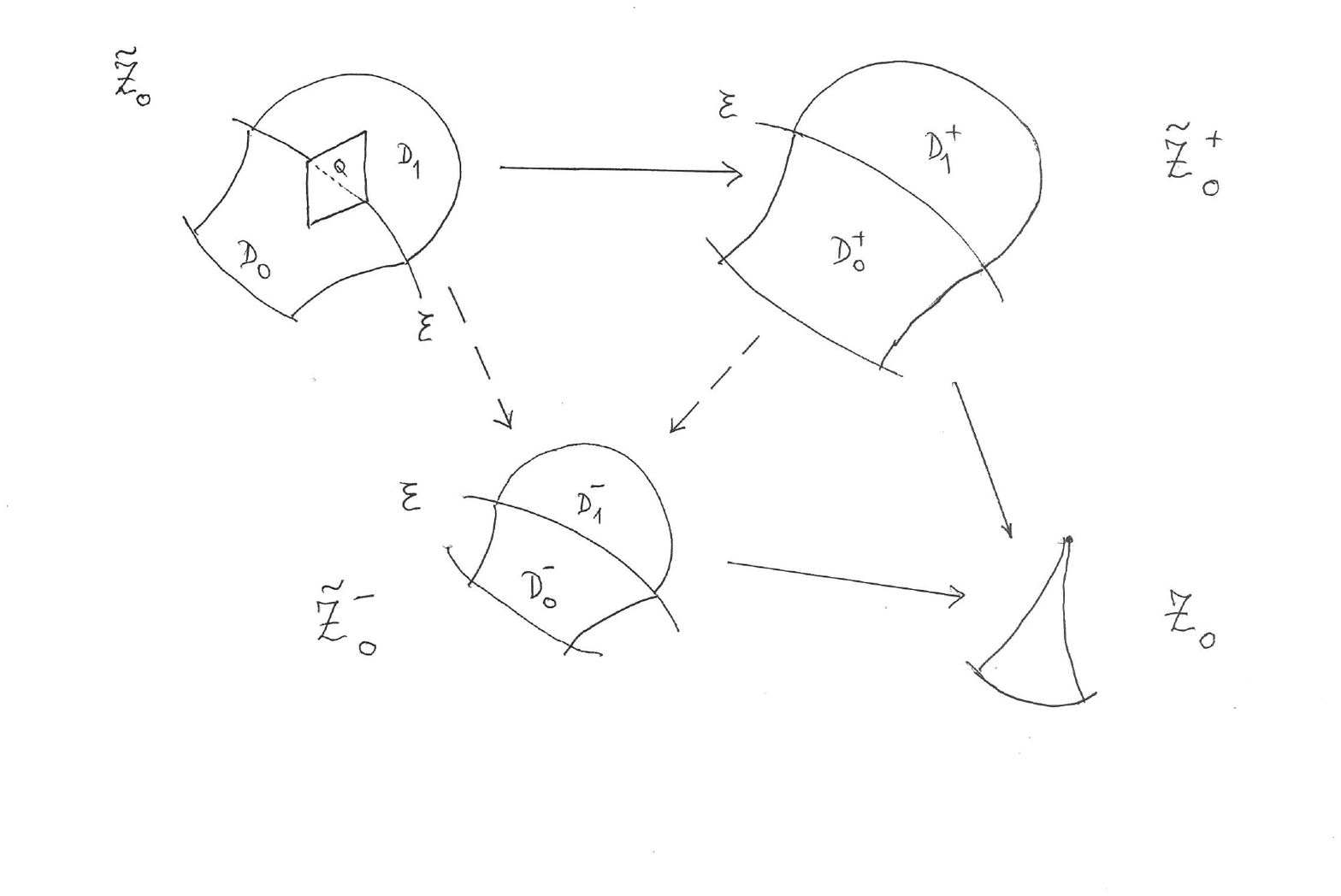}}
\caption{The fibres $\tsZ^{\sharp}$ and $\sZ_0$}
\label{fig:fibres Z etc.}
\end{figure*}

The crucial part of the proof is Theorem \ref{sst},
where we exhibit a surface $D_1$ in $\tsG_{\sE}$
whose transform under a natural flop is a weak del Pezzo
surface. It follows (this is Corollary \ref{HS redux}) that
that the fibre
of $\sZ\to[\hsY/W]$ over the origin of the affine space $[\hsY/W]$ 
has a simply elliptic singularity.
We then recover the result of Helmke and Slodowy,
that $\sZ\to[\hsY/W]$ is a miniversal deformation of this singularity.
In some sense, this is a ``top down'' approach which depends upon being
able to recognize a singularity from its resolution.

\begin{remark} In particular, we have made a non-trivial
birational modification of $\tsG_{\sE}$,
centred on a locus in the boundary $\tsG_{\sE}-\tsG^{ss}_{\sE}$,
which is of a kind that is dictated by the needs
of birational geometry but which has no obvious modular interpretation,
in order to arrive at $\tsZ^-$. 

However, this is a common phenomenon.

For example, consider the coarse moduli space $A_g$ of principally
polarized abelian varieties over $\C$ and its perfect compactification $A_g^P$.
Then [SB16] the exceptional locus of the $\Q$-divisor $12M-D$, where $M$
is the bundle of modular forms of weight $1$ and $D$ is the boundary,
is the image of $A_1^P\times A_{g-1}^P$ in $A_g^P$, the locus of
ppav's with an elliptic factor. So, if $7\le g\le 10$, then $A_g^P$
has terminal singularities [AS16] and is of general type, while
$K_{A_g^P}\sim (g+1)M-D$ is not nef. Therefore
the first step in running the Minimal Model Program on $A_g^P$
is to contract and flip the extremal ray generated by the curve $A_1^P\times\{x\}$,
where $x$ is a point $A_{g-1}^P$. This is a non-trivial birational modification
which is dictated by the MMP
and which is centred in a locus that penetrates the interior $A_g$
of $A_g^P$ that has no obvious modular interpretation, unless the
notion of ppav can be modified.
\end{remark}

It is important to point out that, beyond the work of Brieskorn and Looijenga
mentioned above, and that of Pinkham [P], what we do builds very much 
on the results of Helmke and Slodowy, and, to an equal extent,
on those of Friedman and Morgan [FM].

We will revisit the results of this paper in [GS], where amongst other
things we will
analyse the
moduli stacks in a neighbourhood of infinity by using the Tate curve,
give uniform descriptions of the combinatorics in
terms of the affine Weyl groups and extend the main theorem of this paper
to the cases where $5 \leq d \leq 9$.

Finally, if $H$ is a group then $\B H$ denotes its classifying stack
and $[\sX]$ denotes the geometric quotient of an algebraic stack
$\sX$ when that quotient exists.
\end{section}
\section*{Acknowledgements} We are very grateful to Dougal Davis for
many conversations on these matters and particularly for pointing out an
error in our proof of the previous version of Proposition \ref{6.4}.
We are also grateful
to Igor Dolgachev for his question about unipotent singularities
and to Simon Salamon for many kindnesses, which include drawing the picture.
\begin{section}{The BGSS construction}
To set the notation
we begin by recalling a little of the 
machinery of reductive groups in a way that emphasizes
the flag variety, and its intersection theory,
of an affine group rather than the group itself.
This approach
is very well known, and details can be found in,
for example, [MOSW].

As the starting point of the discussion,
fix a smooth connected 
reductive\footnote{This assumption was omitted in a previous
version of this paper. Otherwise, as was pointed out
by a referee, the field $k$
must be assumed to be perfect.} 
algebraic group $G$ over a field $k$. Then 
[SGA3 XXII, Th. 5.8.1 and Cor. 5.8.3, pp. 228-230]
there is an associated \emph{flag variety} $F$ over $k$ that parametrizes 
the maximal soluble subalgebras (the Borel subalgebras) 
of $\mathfrak g=\Lie G$; it is constructed as a subscheme of the Grassmannian
$Gr(\mathfrak g)$. It is a smooth projective $k$-variety
and is homogeneous under the obvious $G$-action; it is maximal
among the set of homogeneous projective $k$-varieties the stabilizers
of whose points are smooth, and the stabilizer of each point of $F$
is a Borel subgroup
(that is, a maximal connected reduced soluble subgroup of $G$). 
If $G$ has a Borel subgroup $B$ defined over $K$ 
then $F=G/B$, but we do not yet assume that such a 
subgroup exists; that is, we do not assume that $F$ has a $k$-point.
By definition, this last condition is equivalent to
$G$ being \emph{quasi-split}.

Locally in the {\'e}tale topology on $\Sp k$ there is an array of projective
homogeneous $G$-varieties, $2^l$ in number where $l$ is the 
semi-simple rank of $G$,
and $G$-morphisms between them;
in particular, there are $G$-maps $F\to F_1,\ldots,F\to F_l$ that are
{\'e}tale $\P^1$-bundles and there are $G$-maps $F\to X_i$
where $X_i$ is {\emph{minimal}}, that is, of Picard number $1$.
We shall assume that all these varieties, and the morphisms in the lattice,
are defined over $k$; this assumption is fulfilled automatically
if the Dynkin diagram attached to $G$ has no symmetries. Indeed,
the Galois group that acts on the array of varieties above is a subgroup
of the symmetry group of the Dynkin diagram. 
If $G$ is quasi-split, this is, by definition,
the same as $G$ being \emph{split}.

For each minimal variety $X_i$ there is a unique line bundle $\varpi_i$
on $X_i\otimes\bark$ that is a positive generator of $\Pic_{X_i\otimes\bark}$; 
we assume that $\varpi_i$ 
is defined over $K$. This condition holds automatically
if $G$ is split. We also let $\varpi_i$ denote the pullback
of $\varpi_i$ to $F$; then $\varpi_1,\ldots,\varpi_l$ are the
{\emph{fundamental dominant weights}} of $G$ 
and form a $\Z$-basis of $\Pic_F$. (So, for groups of type $E_8$ the bundles 
$\varpi_i$ are always defined over $k$, because here the root lattice 
equals the weight lattice.) They also span the nef cone
of $F$, which is also the effective cone of $F$. 
The semi-simple rank of $G$ is the Picard number of $F$.

Set $\alpha_j=T_{F/F_j}$, the relative tangent bundle.
This is a $G$-linearized line bundle, so defines a class in $\Pic^G_F$,
the group of classes of $G$-linearized line bundles on $F$.
(The rank of this group is the rank of $G$.)
These classes are the {\emph{simple roots}} of $G$.
The {\emph{simple coroots}} are 
$\alpha_1^\vee,\ldots,\alpha_l^\vee$ 
where $\alpha_i^\vee$ is the class of a fibre of $F\to F_i$;
such a curve is a conic (a smooth curve of genus zero).

So some of the various intersection numbers between these curves 
and divisor classes
are given by $(\alpha_i^\vee,\varpi_j)=\delta_{ij}$ and 
$(\alpha_i^\vee,\alpha_i)=2$.

The numbers $1,...,l$ are the nodes of the 
associated Coxeter--Dynkin diagram $D$;
the edges of $D$ are labelled in a way that is
determined by the numbers $(\a_i^\vee,\a_j)$.
This diagram leads to a Coxeter system $(W,S=\{s_1,...,s_l\})$ 
acting on $\Pic^GF$
in the usual way: $s_i(x)=x+(x.\a_i^\vee)\a_i$.
There is also a root datum $(M,M^\vee,\Delta,\Delta^\vee)$
derived from this picture: $M=\Pic^G_F$, $\Delta$ is the $W$-orbit
of the set of simple roots and $\Delta^\vee$ is the $W$-orbit of
the set of simple coroots. More intrinsically, the roots
arise as the $G$-linearized line bundles in a natural
filtration of the tangent bundle $T_F$
and the coroots appear as the classes of curves in the Bruhat
decomposition of $F$.
Conversely, each root datum leads to a unique
split reductive group defined over $\Z$.
(A smooth affine group $G$ over a base $S$ is reductive
if all its geometric fibres are reductive.)

Now assume that $G$ is quasi-split, so split.
Fix a $k$-point of $F$ and rigidify all line bundles
and $G$-linearized line bundles on $F$ at this point.
Define the \emph{weight lattice} $P$ of $G$ by
$P=\Pic^G_F$.
Then there are universal torsors $\sT\to F$
and $\sT^G\to F$ under tori $T$ and $T^G$, respectively,
whose character groups are $\X^*(T)=\Pic_F$ and
$\X^*(T^G)=P$.
Since $G$ is reductive it acts effectively
on $\sT^G$.

If the forgetful homomorphism
$\Pic^G_F\to\Pic_F$ is injective then $G$ is {\emph{semi-simple}}
and if it is an isomorphism then $G$ is {\emph{simply connected}}.
In this case the maps $F\to F_i$ are Zariski $\P^1$-bundles.
A semi-simple group is \emph{simple} if its Dynkin diagram is connected.

\begin{remark}\label{explanation}
It is clear from this that in fact a reductive group does not give rise only
to a root datum, but rather to a \emph{pinned root datum}
[SGA3 XXIII, 1.5]. On the other hand, a del Pezzo surface
$S$ gives rise to a root datum
$(M,M^*,R,R^*)$ in three ways: $M$ is one of 
$K_S^\perp$, $\NS(S)/\Z.[K_S]$ and $\NS(S)$ 
and $R=R^*$ is the set of $(-2)$-vectors
in $\NS(S)$ or the image of that set modulo $\Z.[K_S]$,
but this datum is not naturally pinned, because of the
existence of monodromy on $\NS(S)$.
This explains why the construction to be given in this paper
leads to del Pezzo surfaces that are \emph{split},
in the sense that this monodromy is trivial.
To put it another way, when the Dynkin diagram
has no symmetry, as in the case of $E_8$, there
is no monodromy on $\NS(F)=\Pic_F$
and whether $F$ has a $k$-point is a non-trivial issue.
However, del Pezzo surfaces of degree $1$ always have a 
naturally defined $k$-point
(the base point of the anticanonical system)
while the monodromy is significant.
\end{remark}

{\bf{For the rest of this paper $G$ will be a split
simply connected simple reductive group.}}

Define the torus $T$ by the formula $\X^*(T)=\Pic^G_F$, or
$T=\Hom(\Pic^G_F,\GG_m)$. By assumption, $T$ is split but
is not necessarily isomorphic to a subgroup of $G$.

Next, we recall the BGSS construction.

Define the incidence subvariety $\tG$ of $G\times F$ by
$$\tG=\{(g,x)\vert g(x)=x\}.$$
Then $pr_2:\tG\to F$ is smooth, since $pr_2^{-1}(x)$ is
just the stabilizer $\Stab_G(x)$, and $pr_1:\tG\to G$
is projective. Moreover, there is a smooth
morphism $\pi:\tG\to T$ given by
$$\pi(g,x)(L)= (\phi:g^*L\to L)(x),$$
where the isomorphism $\phi:g^*L\to L$
is part of the data of the $G$-linearization of the line bundle $L$.

\begin{remark} Fix a $k$-point
$x$ on $F$. Take $B=pr_2^{-1}(x)=\Stab_G(x)$. Then $\pi$
restricts to a surjective homomorphism $B\to T$. This homomorphism
is split, and choosing a splitting of it gives the standard set-up
$T\inj B\inj G$ of the theory of reductive groups.
The torus $T$ that we have taken is an \emph{abstract Cartan subgroup}.
\end{remark}

In fact, $(pr_2,\pi):\tG\to F\times T$ is smooth; each fibre
is a translate of the unipotent radical of a Borel subgroup of $G$.

Let $\tG\to X\to G$ be the Stein factorization of $pr_1$.
Then $X=\Sp\G(\tG,\sO_{\tG})$ is the affine hull of $\tG$,
so that $\pi:\tG\to T$ factors through $X$.

From the definitions of $W$ and $T$ there is an action of $W$ 
on $T$ and there is a commutative square
$$\xymatrix{
{\tG}\ar[r]^{pr_1}\ar[d]_{\pi} & {G}\ar[d]^{\rho}\\
{T}\ar[r]&{[T/W]}
}$$
The key point to be proved
here is that the composite morphism $\tG\to T\to [T/W]$ factors through 
$pr_1:\tG\to G$; this can be verified
over $\bark$, where it follows from (1) the existence of a maximal torus
$T_1$ in $G$ that is isomorphic to $T$ and (2) the fact
that the pair $(T_1,N_G(T_1))$ is a slice
to the adjoint action of $G$ on itself.

\begin{proposition}\label{symplectic} The morphism $\pi:\tG\to T$
is relatively symplectic and the canonical class $K_{\tG}$
is trivial.
\begin{proof}
The smooth morphism $pr_2:\tG\to F$ gives, at each point $y$ of $\tG$,
an exact sequence
$$0\to \mathfrak b\to T_{\tG}(y)\to\mathfrak g/\mathfrak b\to 0.$$
So there is an exact sequence
$$0\to\mathfrak u\to T_\pi(y)\to \mathfrak g/\mathfrak b\to 0,$$
where $\mathfrak u$ is the Lie algebra of the unipotent radical 
$U$ of $B$.
The Killing form shows that this sequence is self-dual,
which gives the relative symplectic structure in question.
Since $K_T\sim 0$, 
Then taking the Pfaffian of
the relative symplectic structure
gives the triviality of the relative canonical class $K_{\tG/T}$;
since $K_T$ is trivial, the proposition is proved.
\end{proof}
\end{proposition}

From the viewpoint of algebraic surfaces, their singularities and
their moduli, the significance of this diagram is that, if the 
edges of the Coxeter--Dynkin diagram
are unlabelled (that is, if the group is of type $ADE$), then,
the localization of the unipotent variety $Uni(G)$ 
at the geometric generic point of the subregular unipotent orbit in $G$
has a du Val singularity of the same type as $G$
and the morphism $\rho:G\to [T/W]$ has normal fibres (this can be proved over $\bark$
as a corollary of Steinberg's theorem that in that context $\rho$ 
has a section). Moreover, except in bad characteristic
(taken to include any prime divisor of $n+1$ in the case of type $A_n$), 
$\rho$ yields a miniversal deformation of this singularity 
and $\pi:\tG\to T$ gives
a miniversal deformation of the minimal resolution of the singularity.
\end{section}
\begin{section}{BGSS for semi-stable bundles 
over uniformizable elliptic curves}
Fix an elliptic curve $E$ over a field $k$ and a simply connected reductive
group $G$ over $k$. The word ``bundle'' will imply
``locally trivial in the {\'e}tale topology''.
A principal $G$-bundle $\Xi\to E$ gives rise to an $F$-bundle
$\sF_\Xi=\sF=\Xi\times^GF$ and line bundles $L_{\varpi}=\Xi\times^G\varpi$
on $\sF$ for any weight $\varpi$.
Conversely, suppose that $p:\sF\to E$ is an $F$-bundle and that
$L_{\varpi_1},\ldots,L_{\varpi_l}$ are line bundles on $\sF$
such that $L_{\varpi_i}$ restricts to $\varpi_i$ on
each geometric fibre of $p$. Then there is a reductive and
simply connected group scheme $G_1$ over $E$
defined as the automorphism group scheme of $\sF$ and all the given
line bundles. If the sheaves $p_*L_{\varpi_i}$ are free on $E$
then $G_1$ is constant and pulls back from a $k$-form of $G$.

Fix $\Xi\to E$. Every section $\sigma$ of $\sF\to E$
defines a co-character $[\sigma]$ by
$([\sigma],\varpi)=\sigma.c_1(L_{\varpi})$.
We say that $[\s]\le 0$ if $[\s].\varpi_i\le 0$ for every
fundamental dominant weight $\varpi_i$ and $[\s]<0$ if
$[\s]\le 0$ and $\s\ne 0$.

\begin{definition}
A $G$-bundle $\Xi\to E$ is semi-stable
if $\ad\Xi$ is semi-stable as a coherent sheaf on $E$
and is unstable if it is not semi-stable.
\end{definition}
This is well known ([HS2] Theorem $1.1.2$) to be equivalent to other definitions.
In particular, $\Xi$ is semi-stable if and only if, for every section
$\sigma$ of $\sF_\Xi\to E$, the associated cocharacter
$[\sigma]$ satisfies $[\sigma]\ge 0$.

Recall also Theorem $1.3.1$ of [HS2]: for any 
unstable principal $G$-bundle $\Xi\to E$,
there is a proper parabolic subgroup $P$ of $G$, a 
Levi subgroup $\Lambda$ of $P$ and a principal
$\Lambda$-bundle $\lambda\to E$ such that $\lambda$ is semi-stable,
$\Xi=\lambda\times^{\Lambda}G\to E$ and $\Lambda$ is maximal among
Levi subgroups with this property. Moreover, $\Lambda$ is unique up to
$G$-conjugacy and $\lambda$ is determined up to isomorphism
of principal bundles.

This is proved as a consequence of a unique decomposition
$\ad\Xi=\ll\oplus\uu\oplus\uu^*$ of vector bundles,
where $\ll=\ad\lambda$ consists of summands of degree $0$,
$\uu$ consists of summands of positive degree and
$\ll\oplus\uu=\ad(\lambda\times^{\Lambda}P)$.
In turn, this follows from Atiyah's classification of
vector bundles on an elliptic curve.

Let $\sG^{ss}_E$ denote the stack of semi-stable $G$-bundles on $E$
and $\tsG^{ss}_E$ the stack whose objects are pairs $(\Xi,\sigma)$
where $\Xi\in\Ob \sG^{ss}_E$ and $\sigma$ is a section of the $F$-bundle
$\sF_\Xi=\Xi\times^GF\to E$ whose associated cocharacter
$[\sigma]$ is zero.

\begin{lemma}\label{tautology} $\tsG^{ss}_E$ is naturally isomorphic to the stack
$\sB^0_E$ of $B$-bundles $\beta\to E$ whose associated $T$-bundle
$\beta\times^BT\to E$ is of degree $0$.
\begin{proof} This follows from the tautology that giving a section
of $\sF_\Xi$ is the same as giving a reduction of $\Xi$ to a $B$-bundle.
\end{proof}
\end{lemma}

In the rest of
this section we assume, until the end of the proof
of Theorem \ref{Springer}, that the elliptic curve $E$ is
\emph{uniformizable}; that is, that we can write
$E=\GG_m/\langle q^{\Z}\rangle$.
For example, this holds if either $k=\C$ or $k$
is an algebraically closed complete valued field,
$\ch k\ne 2$ and $j(E)$ is not an integer of $k$
[BGR, 9.7]. If $\ch k=2$ then there do exist uniformizable
curves, derived from the Tate curve. 
We will return to the Tate curve in the future.

Fix a co-ordinate $z$ on $\GG_m$; then we claim that there is a morphism
$\tf:\tG\to\tsG^{ss}_E$ defined by
$$\tf(g,x)=((G\times\GG_m)/\langle\gamma\rangle,\ \{x\}\times E)=(\Xi,\sigma),$$
where $\gamma(h,z)=(ghg^{-1},qz)$. 
The $F$-bundle $\sF_\Xi$ associated to the $G$-bundle $\Xi=\tf(g,x)$
is the quotient $(F\times\GG_m)/\langle\gamma\rangle$
where $\g(y,z)=(g(y),qz)$ for any point $y\in F$,
so that indeed the image of $\{x\}\times E$ is a section of $\sF_\Xi$.

To establish the claim that $\tf(g,x)$ is a point of $\tsG^{ss}_E$,
note that the cocharacter $[\sigma]$ of 
the section $\sigma=\{x\}\times E$ of $\sF_\Xi$
vanishes by continuity, because after specializing to 
the identity element $g=1$ of $G$ it is certainly zero.

\begin{lemma} There is a commutative square
$$\xymatrix{
{\tG}\ar[r]^{\tf}\ar[d]_{\pi} & {\tsG^{ss}_E}\ar[d]^{\tpi}\\
{T}\ar[r]^{\phi}&{Y}
}$$
where $\phi$ is the unramified covering induced from the
uniformization $\GG_m\to E\cong\Pic^0(E)$.
\begin{proof} This follows immediately from the definitions of the 
morphisms involved.
\end{proof}
\end{lemma}

\begin{lemma} $\tpi$ is smooth.
\begin{proof} By Lemma \ref{tautology} we can identify
$\tsG^{ss}_E$ with the stack $\sB^0_E$, and then $\tpi$
is the morphism that maps a $B$-bundle $\beta\to E$ 
to its associated $T$-bundle $\tau=\beta\times^BT\to E$.
The smoothness of $\tpi$ follows from the
surjectivity of the map $H^1(E,\ad\beta)\to H^1(E,\ad\tau)$
of tangent spaces,
which follows in turn from 
the surjectivity of the morphism $\ad\beta\to\ad\tau$
of coherent sheaves on $E$ and the fact that
$E$ is $1$-dimensional.
\end{proof}
\end{lemma}

\begin{lemma}
$[T/W]\to [Y/W]$ is {\'e}tale in an analytic neighbourhood of $0_Y$
and $[Y/W]$ is smooth in a Zariski neighbourhood of $0_Y$.
\begin{proof}
The morphism $\phi:T\to Y$ is the quotient by the group $P$,
which acts freely on $T$,
so $\phi$ is {\'e}tale. Moreover, the group $\tW=P\rtimes W$
acts on $T$ and the stabilizer $\Stab_{\tW}(1_T)$,
where $1_T$ is the identity point of $T$, is $W$,
which equals $\Stab_W(0_Y)$.
The first part of the lemma is now established,
and the rest follows from the smoothness of $[T/W]$.
\end{proof}
\end{lemma}

Now consider the commutative diagram
(``the basic comparison'')
$$\xymatrix{
{\tG}\ar[r]^{\tf}\ar[d]_{pr_1} & {\tsG^{ss}_E}\ar[d]^{s}\\
{G}\ar[r]^{f}& {\sG^{ss}_E}
}$$
where $f(g)=(G\times\GG_m)/\langle\g\rangle$ with
$\g(h,z)=(ghg^{-1},qz)$.
We shall show, in Theorem \ref{classical} below, that this is Cartesian
when restricted to the analytic neighbourhood $\sU$ of the identity 
element $e$ of $G$ that is described in the next Lemma.

\begin{lemma} There is an analytic (classical or rigid) 
open neighbourhood $\sU$ of
$e$ in $G$ on which $f$ is {\'e}tale, and so surjective.
\begin{proof} 
The codifferential $f^*$ is, via Serre duality, a map
$H^0(E,(\ad\Xi)^*)\to\mathfrak g^*$.
For the trivial $G$-bundle 
$\Xi=E\times G$
this is, by inspection, an isomorphism.
So $f$ is smooth at $e$, and the result follows.
\end{proof}
\end{lemma}

Note that $\sU$ maps to an analytic (classical or rigid)
neighbourhood of the trivial
bundle $E\times G$ in $\sG^{ss}_E$.

Recall that $pr_1$ and $s$ are projective, generically finite
and dominant. In fact, $s$ is finite over the locus $\tsG^{ss}_{E,reg}$
of regular semi-stable bundles (those whose automorphism group is of
minimal dimension, namely the rank of $G$).

\begin{lemma} $\deg s =\# W =\deg pr_1$.
\begin{proof} We can assume that there are subgroups $T\inj B\inj G$ as in
the usual set-up for split reductive groups, with a surjection $B\onto T$. 
Then there is a commutative diagram
$$\xymatrix{
{Y}\ar[r]\ar[dr]_{=}&{\tsG^{ss}_{E,reg}}\ar[r]\ar[d]&{\sG^{ss}_{E,reg}}\\
& {Y}
}$$
arising from the identification made above of $\tsG^{ss}_E$ with the
stack $\sB^0_E$. 
Now $\sG^{ss}_{E,reg}$ has a geometric quotient that can be identified with
$[Y/W]$ and the fibres of $\tsG^{ss}_{E,reg}\to Y$ are
point over the regular locus in $Y$. That is, over the maximal
open subvariety $Y^0$ of$Y$ on which $Y$ acts freely.

In other words, the square
$$\xymatrix{
{\tsG^{ss}_{E,reg}}\ar[r]\ar[d]&{\sG^{ss}_{E,reg}}\ar[d]\\
{Y}\ar[r]&{[Y/W]}
}$$
is Cartesian when restricted to $Y^0$
and this is enough.
\end{proof}
\end{lemma}

\begin{theorem}\label{classical} 
The square
$$\xymatrix{
{\tG}\ar[r]^{\tf}\ar[d]_{pr_1} & {\tsG^{ss}_E}\ar[d]^s\\
{G}\ar[r]^{f}& {\sG^{ss}_E}
}$$ in the ``basic comparison'' is Cartesian when restricted to 
the neighbourhood $\sU$ of $e$ in $G$.
\begin{proof}
Consider the diagram
$$\xymatrix{
{\tG}\ar[drr]^{\tf}\ar[dr]^r\ar[ddr]_{pr_1}&&\\
& {Z}\ar[r]_{h'}\ar[d]^q & {\tsG^{ss}_E}\ar[d]^s\\
& {G}\ar[r]^f & {\sG^{ss}_E}
}$$
where $Z$ is the fibre product. Since $\deg pr_1=\deg s$, it
follows that $r$ is birational; clearly, $r$ is proper.
Also, $h'$ is smooth on $\sU_1=q^{-1}(\sU)$, since $f$ is smooth
on $\sU$,
so that $\sU_1$ is smooth. Set $\sU_2=r^{-1}(\sU_1)$,
an  open subvariety of $\tG$.
So $r:\sU_2\to\sU_1$ is a proper birational morphism
of smooth analytic $k$-varieties.

By Proposition \ref{symplectic}
the canonical divisor class $K_{\tG}$ is trivial,
and then $K_{\sU_2}\sim 0$.

\begin{lemma} If $g:X\to Y$ is a proper birational morphism of 
smooth (rigid analytic or classical) $k$-varieties and
$K_X\sim 0$, then $g$ is an isomorphism.
\begin{proof} This is very well known in characteristic zero
but possibly less so in positive characteristic.
So we give a proof.

There is an open subvariety $V$ of $Y$ whose complement is
of codimension at least $2$ over which $g$ is an isomorphism.
Therefore the trace map $tr:g_*\omega_X\to \omega_Y$ is an 
isomorphism on $V$; since $g_*\omega_X\cong\sO_Y$
it follows that $tr$ is an isomorphism and that
$\omega_Y\cong\sO_Y$. Then the natural homomorphism
$g^*\omega_Y\to\omega_X$ is also an isomorphism.
Therefore the determinant of the derivative $dg:T_X\to g^*T_Y$
is an isomorphism, and so $dg$ is an isomorphism.
Therefore $g$ is smooth. Since it is also proper and birational
it is an isomorphism.
\end{proof}
\end{lemma}

In particular, the restriction of $r$
to $\sU_2$ is an isomorphism and we are done.
\end{proof}
\end{theorem}

\begin{corollary}\label{sections exist} 
$s:\tsG^{ss}_E\to \sG^{ss}_E$ is proper and surjective.
\begin{proof}
By Theorem \ref{classical}
the morphism $s$ is surjective in a neighbourhood of the
trivial bundle in $\sG^{ss}$. The properness
follows from the valuative criterion and the
fact that, if $\{\Xi_t\}$ is a family of semi-stable
$G$-bundles over $E$, then a section $\s_t$ of 
$\sF_{\Xi_t}$
whose cocharacter $[\s_t]$ vanishes
can only specialize to a section $[\s_0]$
of $\sF_{\Xi_0}$ such that $[\s_0]\le 0$.
Since $\Xi_0$ is semi-stable, $[\s_0]=0$.
\end{proof}
\end{corollary}

\begin{theorem}\label{Springer}
In a neighbourhood of the origin in $[Y/W]$
the squares
$$\xymatrix{
{\tsG^{ss}_E}\ar[r]^s\ar[d] & {\sG^{ss}_E}\ar[d]\\
{Y} \ar[r] & {[Y/W]}
}
{\ and\ \ }
\xymatrix{
{\tG}\ar[r]\ar[d] & {G}\ar[d]\\
{T}\ar[r] & {[T/W]}
}$$
are smoothly equivalent.
\begin{proof}
Define $\sH=\sG^{ss}_E\times_{[Y/W]}[T/W]$ and
$\tsH=\tsG^{ss}_E\times_YT$, so that
the square
$$\xymatrix{
{\tsH}\ar[r]\ar[d] & {\sH}\ar[d]\\
{T}\ar[r] & {[T/W]}
}$$
is the pull-back under ${[T/W]}\to [Y/W]$
of the square
$$\xymatrix{
{\tsG^{ss}_E}\ar[r]^s\ar[d] & {\sG^{ss}_E}\ar[d]\\
{Y} \ar[r] & {[Y/W].}
}$$

We have shown that, near $e$, the square
$$\xymatrix{
{\tG}\ar[r]\ar[d] & {G}\ar[d]\\
{\tsG^{ss}_E}\ar[r] & {\sG^{ss}_E}
}$$
is Cartesian; the same argument shows that, near $e$, $G\to\sH$
is smooth and that
$$\xymatrix{
{\tG}\ar[r]\ar[d] & {G}\ar[d]\\
{\tsH}\ar[r] & {\sH}
}$$
is Cartesian in a neighbourhood of $e$. So,
in a neighbourhood of $e$,
$$\xymatrix{
{\tG}\ar[r]\ar[d] & {G}\ar[d]\\
{T}\ar[r] & {[T/W]}
}$$
is smoothly equivalent to
$$\xymatrix{
{\tsH}\ar[r]\ar[d] & {\sH}\ar[d]\\
{T}\ar[r] & {[T/W].}
}$$
Comparing these two descriptions of
the last square gives the result.
\end{proof}
\end{theorem}

We shall see later that this result 
fails when $E$ is supersingular in characteristics
$2,3$ and $5$ and $G=E_8$.

We can now tie up some loose ends concerning semi-stability.
We no longer assume the curve $E$ to be
uniformizable.

We have recovered, as Corollary \ref{sections exist},
the well known fact that the first projection
$s:\tsG^{ss}_E\to\sG^{ss}_E$ is surjective and proper when the base is $\Sp\C$.
This is true over any base.

\begin{lemma} The projection
$s:\tsG^{ss}_{\sE}\to\sG^{ss}_{\sE}$ is surjective and projective.
\begin{proof} It follows from Corollary \ref{sections exist}
that $s:\tsG^{ss}_{E}\to\sG^{ss}_{E}$ is surjective 
if $E$ is defined over $\C$, and so over
any field of characteristic zero.
Suppose then that $\Xi$ is a semi-stable $G$ bundle over $E$
in characteristic $p$; lift $E$ and $\Xi$
to ${\widetilde{\Xi}}\to\tE$ in
characteristic zero. There is then a section
$\ts$ of $\sF_{\widetilde{\Xi}}$ with $[\ts]=0$.
Specializing back to characteristic $p$
gives a section $\s$ of $\sF_\Xi$ in characteristic $p$
with $[\s]\le 0$; since $\Xi$ is semi-stable,
$[\s]=0$.

So $s$ is surjective. Properness follows from the valuative
criterion in the same way and then projectivity is a 
consequence of the projectivity of the Hilbert scheme.
\end{proof}
\end{lemma}

\begin{proposition}\label{semi-stable vs. unstable}
If the base is an algebraically closed field,
then $\Xi$ is semi-stable if and only if $\sF_\Xi$ has
a section $\sigma$ with $[\sigma]=0$
and $\Xi$ is unstable if and only if $\sF_\Xi$
has a section $\s$ with $[\s]<0$.
\begin{proof} Suppose that $\Xi$ is semi-stable. Then the existence
of a suitable $\sigma$ is the surjectivity to which we have just alluded.
Conversely, the existence of a section $\sigma$ such that $[\sigma]=0$
leads to a description of $\ad\Xi$ as an extension of line bundles
of degree $0$, which implies its semi-stability.

According to [HS3], 3.3, every unstable
$G$-bundle $\Xi_0$ can be deformed to a semi-stable
bundle $\Xi_t$. Then $\sF_{\Xi_t}$ has a section $\s_t$
with $[\s_t]=0$; specializing to $t=0$ gives
a section $\s_0$ of $\sF_{\Xi_0}$ with $[\s_0]\le 0$.
Since $\Xi_0$ is unstable, $[\s_0]\ne 0$.
\end{proof}
\end{proposition}
\begin{remark}\label{explain}
There is a special case of
Theorem 1.3.1 of [HS2] that we can now spell out.

Suppose that $\Xi$ can be reduced to a semi-table $\Lambda$-bundle
where $\Lambda$ is a Levi subgroup of a maximal parabolic subgroup
$P$ of $G$. Let $\Xi_P=\Xi/P=\Xi\times^G G/P\to E$
be the associated $(G/P)$-bundle. Then $\Xi_P\to E$ has a section
$\s$ whose associated cocharacter $[\s]$ 
is a negative integer multiple of $\a^\vee$, where $\a^\vee$ is the simple
coroot belonging to $P$.
\end{remark}
\end{section}
\begin{section}{Some relative compactifications of $\tsG^{ss}_{\sE}$}
\label{compactification}
We denote by $\sT^0$ the stack of $T$-bundles over $\sE$ of multi-degree zero,
so that $\sT^0\cong \sY\times \B T$, where $\sY=\Hom(P,\sE)$
and, for a group $H$, $\B H=\{*\}/H$ is the classifying stack of $H$-bundles.
There are morphisms $s:\tsG^{ss}_{\sE}\to\sG^{ss}_{\sE}$ 
and $\tpi:\tsG^{ss}_{\sE}\to \sY$ given by
$\tpi(\Xi,\sigma)(\varpi)=L_{\varpi}\vert_{\sigma}$
and these stacks and morphisms fit into a commutative square,
analogous to the square in the BGSS construction,
$$\xymatrix{
{\tsG^{ss}_{\sE}}\ar[r]^{s}\ar[d]_{\tpi} & {\sG^{ss}_{\sE}}\ar[d]\\
{\sY}\ar[r] & {[\sY/W]}
}$$
where $s=pr_1$ is projective (in particular, representable),
by Proposition \ref{6.1} below,
and $[\sY/W]$ is the geometric quotient,
relative to $\sM_{Ell}$, of the open substack $\sG^{ss}_{\sE,reg}$
of {\emph{regular}} semi-stable bundles. Here ``regular'' means that
the automorphism group has minimal dimension $l$, the rank of $G$.

It is well known that the classifying morphism
$\sG^{ss}_{\sE,reg}\to [\sY/W]$ extends to a
morphism $\g:\sG^{ss}_{\sE}\to [\sY/W]$
given by sending a bundle to its $S$-equivalence class.

We now recall various relative compactifications of $\tsG^{ss}_{\sE}$.
There are several that are relevant, but the most useful for us here,
because of its smoothness properties,
will be denoted by $\tsG_{\sE}$ or $\tsG^{KM}_{\sE}$; it
is based on the stack of stable maps introduced by Kontsevich
(although it is also appropriate to attach the name of Mori).
That is, it relies on enlarging the source of a map, rather than on Drinfel'd's
idea, which we now recall, of enlarging the target.

Define $R=\oplus_\varpi H^0(F,\varpi)$,
the Cox ring of $F$, the sum being taken over all dominant weights.
Let $T$ denote the torus introduced previously
and consider the singular and non-separated
stack $\bF=(\Sp R)/T$.
There is an open $G$-equivariant embedding $F\inj\bF$
and Drinfel'd's idea, which is described in Section $1$ of [BG], and recalled
in more detail below, is to embed $\sF_\Xi$ into the $\bF$-bundle
$\Xi\times^G\bF\to E$ and consider sections of this.
Of course, any projective
homogeneous $G$-variety $X$ has a similar enlargement $X\inj\bX$.
This Kontsevich-Mori compactification turns out to be 
slightly wrong for our purposes 
(roughly speaking, it needs to be contracted and then flopped)
but
the extra information that it contains turns out to be
crucial for the proof of our main result. The Drinfel'd compactification
$\tsG^D_{\sE}$, on the other hand, 
is too small and too singular.

Let $\sC_{pre}\to\Sp\Z$ denote the stack whose objects over 
a point $S$ are $1$-marked
pre-stable curves $C\to S$ whose canonical model is an 
elliptic curve over $S$.
So, if $S$ is a geometric point, then $C$ is reduced with normal crossings,
its dual graph is a tree, every irreducible component except one,
say $C_1$, a a copy of $\P^1$ and $C_1$ is an elliptic curve. 
We refer to these as ``elliptic curves with rational tails''
and to $\sC_{pre}$ as the stack of elliptic curves with tails. Note that
$\sC_{pre}$ is a smooth stack over $\Z$ and that 
its discriminant, the locus of singular curves,
is a divisor $\Delta=\Delta(\sC_{pre})$ in $\sC_{pre}$ with normal crossings.
The complement $\sC_{pre}-\Delta$ is naturally isomorphic to
$\sM_{Ell}$ and there is a retraction $\sC_{pre}\to\sM_{Ell}$
given by sending a pre-stable curve to its canonical model.
This retraction is smooth, since it is surjective on tangent spaces.
 
Let $\tsG_{\sE}=\tsG^{KM}_{\sE}$ 
be the stack whose objects over a scheme $S$ are triples
$(\Xi,C,\sigma:C\to\Xi\times^GF)$, where $\Xi$ is a
$G$-bundle over $E=\sE\times_{\sM_{Ell}} S$, 
$C$ is an $S$-object of $\sC_{pre}$ and $\sigma$
is a stable map (in the sense of Kontsevich) such that the
composite $C\to\Xi\times^GF\to E$ is the contraction of $C$ to its
canonical model and $\deg\sigma^*L_{\varpi_i}=0$
for each fundamental dominant weight $\varpi_i$.
(Here, degree means total degree, the sum of the degrees
on each component of $C$.)
Let $s:\tsG\to\sG$ be the forgetful map $(\Xi,C,\sigma:C\to\Xi\times^GF)\mapsto \Xi$.

\begin{proposition}\label{*}\label{6.1}
\part[i] $s:\tsG_{\sE}\to\sG_{\sE}$ is proper and has finite 
relative automorphism group schemes.

\part[ii] $pr_2:\tsG_{\sE}\to\sC_{pre}$ is smooth.

\part[iii] There is a smooth morphism $\pi:\tsG_{\sE}\to \sY$
that extends the morphism $\tpi:\tsG^{ss}_{\sE}\to \sY$ described previously.

\part[iv]\label{smoothness} 
$(pr_2,\pi):\tsG_{\sE}\to\sC_{pre}\times_{\sM_{Ell}}\sY$
is also smooth.

\part[v]\label{representable} $s$ is representable when restricted
to the locus in $\tsG_{\sE}$ where, on each geometric component of $C$,
either $\sigma$ is of degree $1$ onto its image or $\sigma$ is constant.
\begin{proof} The only things which are neither obvious nor well known
are \DHrefpart{iii} and \DHrefpart{iv}. 

For \DHrefpart{iii}, it is enough to consider the universal
curve $\G\to\sC_{pre}$ and then notice that the contraction
$\G\to \sE\times_{\sM_{Ell}}\sC_{pre}$ is a projective and birational morphism of
smooth stacks. It is well known that for such a morphism
there is a blowing-down morphism from $\GG_m$-bundles on $\G$
to $\GG_m$-bundles on $\sE\times_{\sM_{Ell}}\sC_{pre}$, given by pushing forward
divisor classes; this extends to a blowing-down morphism
from $T$-bundles on $\G$
to $T$-bundles on $\sE\times_{\sM_{Ell}}\sC_{pre}$ for any split torus $T$,
by induction on the rank of $T$.
So $\pi$ exists; its smoothness is an immediate consequence
of the fact that curves are $1$-dimensional.

\DHrefpart{iv}: Suppose that $k$ is an algebraically closed field
and that $C$ is a $k$-point of $\sC_{pre}$ whose canonical model
is the elliptic curve $E$. Since $H^1(E,\sO_E)\to H^1(C,\sO_C)$
is an isomorphism, it is enough to show that
the morphism $\sB_C\to\sT_C$ from the stack of $B$-bundles on $C$ to 
the stack of $T$-bundles
on $C$ is smooth.
As usual, the obstruction to smoothness lies in a group of the 
form $H^2(C,\ad\a)$,
where $\a\to C$ is a principal bundle under a unipotent group,
so vanishes.
\end{proof}
\end{proposition}

We shall see that \ref{representable} is enough to show that $s$ is representable
over some neighbourhood of the locus of regular or subregular unstable $G$-bundles.

The objects of Drinfeld's compactification $\tsG^D_{\sE}$
are described as follows; we refer to [BG, section 1] for the details
of what is described in the next paragraph.

Fix an elliptic curve $E$ and
a $G$-bundle $\Xi\to E$, with
$q:\sF=\sF_\Xi=\Xi\times^GF\to E$ the associated $F$-bundle.
For every dominant weight $\varpi$ there is a line bundle
$L_\varpi$ on $\sF$ and vector bundle $V_\varpi=q_*L_\varpi$ on $E$.
A reduction of $\Xi$ to a $B$-bundle is a line sub-bundle
$M_\varpi$ of $V_\varpi$ (the line generated by a vector of highest weight)
such that the set of all subsheaves
$M_\varpi\inj V_\varpi$, as $\varpi$ ranges over all dominant weights,
satisfies the Pl{\"u}cker relations. An object of $\tsG^D_{\sE}$
consists of a $G$-bundle $\Xi$ over an elliptic curve $E$
and a collection of subsheaves
$\{M_\varpi\inj V_\varpi\}_\varpi$ where $M_\varpi$ is invertible, but not
necessarily a sub-bundle,
that satisfies the Pl{\"u}cker relations; we also demand that
the associated $T$-bundle should have degree $0$.
(This $T$-bundle is constructed as follows:
objects of $\tsG^D_{\sE}$ are identified with
sections of the $\bF$-bundle $\bsF=\Xi\times^G\bF\to E$ that meet the open
subscheme $\sF$ of $\bsF$. The open embedding
$F\inj\bF$ induces an isomorphism
$\Pic^G_{\bF}\to \Pic^G_F$, so there is a natural $T$-bundle associated
to $\bsF$; this is to be of degree $0$.)

The projection $\tsG^D_{\sE}\to\sG_{\sE}$ is projective,
while $\tsG^{KM}\to\sG$ is also proper but can have
non-trivial, but finite, relative automorphism groups.

\begin{proposition}\label{natural} There is a natural morphism 
$\tsG_{\sE}\to\tsG^D_{\sE}$
relative to $\sG_{\sE}$.
\begin{proof} An object of $\tsG_{\sE}$ gives sub-line bundles $M'_\varpi$
of $V_\varpi$ on a pre-stable curve $C$ where each $M'_\varpi$
has degree $0$; pushing these sheaves
$M_\varpi$ forward to $E$ gives an object of $\tsG^D_{\sE}$.
\end{proof}
\end{proposition}
\end{section}
\begin{section}{Some deformations}
and the stack $\tsG^+_{\sE}$
\label{some charts}\label{G+}\label{previous}\label{6}
Assume that $G$ is a split and simply connected reductive group 
of type $E_l$,
where $l=5,6,7$ or $8$, and that $T\subset G$ is a copy of
$\Hom(\Pic^G_F,\GG_m)$, a maximal torus in $G$.
Fix also a Borel subgroup $B$ containing $T$. 
We shall assume, as we may, that these data are all given over
$\Sp\Z$. For $E_6,E_7$ and $E_8$
we shall number the nodes of Dynkin diagrams as in Bourkaki's
planches: the branch node is numbered $4$
and the node adjacent to it on the long arm (or one of the two longer arms
in the case of $E_6$) is numbered $5$. For $E_5=D_5$ we also number
the branch node by $4$, but $5$ will refer to a node on one of the short
arms. Here is the diagram for $E_8$.
$$\xymatrix{
&&{\bullet}^{2}\ar@{-}[d]&&&&\\
{\bullet}_1\ar@{-}[r]&{\bullet}_3\ar@{-}[r]&{\bullet}_4\ar@{-}[r]&{\bullet}_5\ar@{-}[r]
&{\bullet}_6\ar@{-}[r]&{\bullet}_7\ar@{-}[r]&{\bullet}_8
}$$
Recall, from [FM], the construction
of miniversal deformation spaces for certain unstable
$G$-bundles over an elliptic curve $E\to S$,
where $S$ is some scheme.
Suppose that $P\supset B$ is a maximal parabolic subgroup,
with a Levi factor $\Lambda=P\cap P^-$, where
$P^-\supset B^-$ are opposite to the pair
$P\supset B$.
Fix a $\Lambda$-bundle $\lambda\to E$
such that

\noindent $(i)$ $\lambda$ is semi-stable,

\noindent $(ii)$ $\Xi=\lambda\times^{\Lambda}G$ is an unstable $G$-bundle and

\noindent $(iii)$ $\lambda$ is regular, in the sense 
that for every geometric point of $S$ the 
automorphism group of $\lambda$ is of
minimal dimension amongst all $\Lambda$-bundles that
satisfy $(i)$ and $(ii)$.

Then consider the stack $\sH$ whose objects are triples $(\lambda,\sP^-,\phi)$,
where $\sP^-$ is a $P^-$-bundle over $E$ and $\phi$ is an isomorphism
$\phi:\sP^-/U^-\to\lambda$, where $U^-$ is the unipotent radical of $P^-$.
The forgetful map $\sH\to\Lambda^{ss}_{\sE,reg}$, where $\Lambda^{ss}_{\sE,reg}$ 
is the stack of regular semi-stable $\Lambda$-bundles over $E$, is represented by an 
affine space bundle over $\Lambda^{ss}_{\sE,reg}$; the fibre over $\lambda$
is the non-abelian cohomology set $H^1(E,{\underline{U}}^-)$, where
${\underline{U}}^-$ is the principal $U^-$-bundle $\lambda\times^{\Lambda} U^-$.
However, it is shown in [FM], section $4$, that this set is naturally an affine space
isomorphic to the cohomology vector space $H^1(E,\uu^-)$, where
$\uu^-$ is the corresponding bundle of Lie algebras. (In {\it{loc. cit.}}
the base is $\Sp\C$ but their argument is valid
over any base. 
So if $E$ is defined over a base scheme $S$
then $H^1(E,{\underline{U}}^-)$ is an affine space bundle over $S$.)
In particular,
$\sH$ is algebraic and smooth.

Consider the morphism $\rho:\sH\to\sG_{\sE}$ given by 
$\rho(\sP^-)=\sP^-\times^{P^-} G$.
Note that there is an action of the centre $\zeta(\Lambda)$ on $\sH$ that
covers the trivial action on $\sG$, coming from the inclusion 
$\zeta(\Lambda)\inj\sP^-$.
The fixed locus of this action is the stack of triples
$(\lambda,\lambda\times^{\Lambda}P^-,\phi_{can})$, which is a copy of $\Lambda^{ss}_{reg}$.
That is, the affine bundle $\sH\to\sL^{ss}_{reg}$ has a $\GG_m$-action
and a section consisting of fixed points for the $\GG_m$-action.

Now suppose that we start with an unstable $G$-bundle $\xi$,
defined over some algebraically closed field $k$, that is
either regular or subregular. Then, according to [HS2], Theorem $5.12$,
$\xi$ is isomorphic, up to translation by a $k$-point of $E$,
to the bundle $\Xi_r=\lambda_r\times^{\Lambda_r}G$ 
where $r=4$ if $\xi$ is regular
and $r=5$ if $\xi$ is subregular, $\lambda_r$ is a
semi-stable $\Lambda_r$-bundle and
$\Lambda_r$ is a Levi subgroup of a maximal parabolic subgroup $P_r$
that is associated to the node numbered $r$
From its description, the group $\Lambda_r$ is defined over $\Z$.

\begin{lemma}
The bundle $\lambda_r$ is defined 
over $\sM_{Ell}$. 
\begin{proof} When $r=4$ this follows from the description
of $\lambda_r$ on pp. 375--376 of [HS2]
and the fact that $\sE$ has a natural line bundle of degree $1$,
corresponding to the origin $0_{\sE}$.

When $r=5$ there is a similar description:
$\lambda_5$ is
determined by a $GL_q \times GL_s$ bundle
$(\eta_q,\eta_s)$ with $\deg\det\eta_s=1$ and $\det \eta_q \cong \det \eta_s^{\otimes 2}$,
where $(q,s) = (5,l-4)$.
\end{proof}
\end{lemma}

Therefore there is a versal deformation space $\szd$ for $\xi$ that is
an affine space bundle, on which $\GG_m$ acts with strictly positive weights, 
over $\Pic^1(\sE)$. The fixed locus of the $\GG_m$-action is a section
and we can projectivize to form $\P(\szd)\to\Pic^1(\sE)$.

Let $\szdss$ be the semistable locus in $\szd$; this is open, and its
complement has codimension $\ge 2$. There is a classifying morphism
$\g:\szdss\to [\sY/W]$.

Consider the induced morphism
$\phi:\szdss\to [\sY/W]\times_{\sM_{Ell}}\Pic^1(\sE)$.
This morphism is constant on
$\GG_m$-orbits and therefore 
factors through a morphism
$\psi:[\szdss/\GG_m]\to [\sY/W]\times_{\sM_{Ell}}\Pic^1(\sE)$.

Note that $[\szdss/\GG_m]$ is an open substack of 
the weighted projective space bundle $\P(\szd)$
over $ [\sY/W]\times_{\sM_{Ell}}\Pic^1(\sE)$
and its
complement is of codimension at least $2$.

For any line bundle $\sL$ on $[\sY/W]\times_{\sM_{Ell}}\Pic^1(\sE)$
that is ample relative to $\Pic^1(\sE)$,
let $\widehat{\sL}\to\Pic^1(\sE)$
be the family of corresponding cones over $[\sY/W]$
obtained by contracting the $0$-section $0_{\sL}$
of the total space of $\sL$. 
So there is a commutative diagram
$$\xymatrix{
{\left(\sL-0_{\sL}\right)}\ar[d]\ar@{^(->}[r]|-{o}&{\widehat{\sL}}\ar[d]\\
{[\sY/W]\times_{\sM_{Ell}}\Pic^1(\sE)}\ar[r]^<<<<<{pr_2}&{\Pic^1(\sE).}
}$$

\begin{lemma}\label{lift} Suppose that $\sL_1$ is a line bundle
on $[\sY/W]\times_{\sM_{Ell}}\Pic^1(\sE)$ that is ample
relative to $\Pic^1(\sE)$.
Then there is a line bundle $\sA$ on $\Pic^1(\sE)$ such that,
if $\sL=\sL_1\otimes pr_2^*\sA$, then
$\phi$ lifts to $\tphi:\szdss\to \sL-0_{\sL}$.
\begin{proof} Write $\sM^0=\psi^*\sL_1$. Then, after
replacing $\sL_1$ by $\sL$ as in the statement of the lemma,
we have $\sM^0=\sO_{\P(\szd)}(n)\vert_{[\szdss/\GG_m]}$
for some $n\ge 1$.

Write $\sM=\sO_{\P(\szd)}(n)$.
Then $\sM^0-0_{\sM^0}=\szdss/\mu_n$
and there is a commutative diagram
$$\xymatrix{
{\szdss}\ar[r]\ar[dr]&{\left(\sM^0-0_{\sM^0}\right)}\ar[d]\ar[r]^{\Psi}&
{\left(\sL-0_{\sL}\right)}\ar[d]\ar@{^(->}[r]|-{o}&{\widehat{\sL}}\ar[d]\\
&{[\szdss/\GG_m]}\ar[r]^>>>>>>>>>{\psi}&
{[\sY/W]\times_{\sM_{Ell}}\Pic^1(\sE)}\ar[r]^<<<<<<{pr_2}&{\Pic^1(\sE).}
}$$
Composing the arrows in the top row
gives the result.
\end{proof}
\end{lemma}

\begin{lemma}\label{affine}
Suppose that $\sT$ is a stack, that $\sX$ and $\sV$ are flat
and affine over $\sT$, that $\sX$ is normal
and that $\sU$ is an open substack of $\sX$
whose complement has codimension at least $2$
in every fibre of $\sX\to\sT$.
Then every morphism $f:\sU\to\sV$ over $\sT$
extends uniquely to a morphism $\sX\to\sV$.
\begin{proof}
Suppose that $T\to\sT$ is a smooth cover and that $T$
is affine. Write $\sU_T=\sU\times_{\sT}T$ etc.
and we have a diagram
$$\xymatrix{
{\sU_T}\ar@{^(->}[r]|-{o}\ar[dr]&{\sX_T}\\
&{\sV_T}
}$$
of schemes, where $\sX_T$ and $\sV_T$ are affine and $\sX_T$ is normal.
Now
\begin{eqnarray*}
Mor_{Sch/T}(\sU_T,\sV_T)&=&Mor_{Rings/\G(\sO_T)}(\G(\sO_{\sV_T}),\G(\sO_{\sU_T}))\\
&=&Mor_{Rings/\G(\sO_T)}(\G(\sO_{\sV_T}),\G(\sO_{\sX_T}))=Mor_{Sch/T}(\sX_T,\sV_T)
\end{eqnarray*}
since $\G(\sO_{\sU_T})=\G(\sO_{\sX_T})$.
So the lemma is proved in the scheme-theoretical
context and then follows in full via descent.
\end{proof}
\end{lemma}

\begin{corollary} The morphism $\g:\szdss\to 
{[\sY/W]\times_{\sM_{Ell}}\Pic^1(\sE)}$
lifts to a morphism $\szd\to{\widehat{\sL}}$.
\begin{proof} This follows from the previous two lemmas.
\end{proof}
\end{corollary}

Looijenga [L1] considered ample $W$-linearized
line bundles $\sL$ on $\sY$
such that the linearization is trivial on the line $\sL(0_{\sY})$
and the polarization defined by $\sL$
is the standard $W$-invariant quadratic form on 
the weight lattice $P$.

He proved that the set of such line bundles is
naturally a torsor under the finite group scheme
$Hom(P/Q,\sE)$ and that
over any complex point of $\sM_{Ell}$
the ring of invariants $\oplus_{n\ge 0} H^0(\sY,\sL^{\otimes n})^W$
is a polynomial ring and is generated by
homogeneous elements whose degrees are the coefficients
of the biggest (highest) root when expressed
in terms of simple roots. (Recall that we are only
considering simply laced groups.)

\begin{proposition}
There exists a line bundle on $\sY$
that satisfies Looijenga's hypotheses.
\begin{proof} Assume that $r\le 8$
and let $\Lambda$ denote the hyperbolic
lattice $\Lambda=\Z e_0\oplus \oplus_1^r\Z e_i$
with $e_0^2=-1, e_i^2=1$ for $i\ge 1$ and 
$e_i.e_j=0$ otherwise. Put $\kappa= 3e_0-\sum_1^re_i$.
Then $P\cong\Lambda/\Z\kappa$
and the exact sequence
$$0\to\Z\kappa\to\Lambda\to P\to 0$$
gives an exact sequence
$$0\to\sY\to Hom(\Lambda,\sE)\to \sE\to 0.$$
For $j=0,1,...,r$ we have $\Z e_j\inj \Lambda$,
so projections $pr_j:Hom(\Lambda,\sE)\to \sE$.
Define $\sA=\sO_{\sE}(0_{\sE})$
and $\sA_j=pr_j^*\sA$. 
Then $\sL=\sA_0^{-1}\otimes\sA_1\otimes\cdots\otimes\sA_r$
is the line bundle required.
\end{proof}
\end{proposition}

That is, each regular or subregular unstable bundle 
$\Xi_4$ or $\Xi_5$ has a versal deformation
space that is an affine bundle over $\Pic^1(\sE)$, 
and after turning off $\Pic^1(\sE)$
the affine space that is the fibre admits a 
classifying morphism to the cone $[\hsY/W]$
corresponding to the line bundle $\sL$.

Now we construct the stack $\tsG^+_{\sE}$.
Fix the $G$-homogeneous spaces $H_4,H_5$ and $H_{4,5}$ 
associated to
the sets of nodes $\{4\}$, $\{5\}$ and $\{4,5\}$,
so that $F\to H_4$ and $F\to H_5$ are $\P^1$-bundles,
$H_4\to H_{4,5}$ and $H_5\to H_{4,5}$ are $\P^2$-bundles
and $F\to H_{4,5}$ is a bundle whose fibre is the flag variety
$SL_3/B$ of type $A_2$.

Given $S\to\sM_{Ell}$,
the $S$-objects of $\tsG^+_{\sE}$ are triples $(\Xi,\s,\tau)$
where $\Xi$ is a $G$-bundle $\Xi\to \sE\times_{\sM_{Ell}} S$, $\s$ a stable map
$\s:C\to \sH_5=\Xi\times^G H_5$ and $\tau$ a section of
$\bsF$ such that the associated $T$-bundle is of degree $0$.
Moreover, we require that under the projection $\bsF\to\bsH_5=\Xi\times^G\bH_5$
the curves $\s$ and $\tau$ should agree over a dense open subset $U$ of 
$\sE\times_{\sM_{Ell}} S$ that meets every geometric fibre of 
$\sE\times_{\sM_{Ell}} S\to S$.
It follows from the definition that $\tsG^+_{\sE}$ is a closed substack
of $\tsG_{\sE}\times_{\sG_{\sE}}\tsG^D_{\sE}$.

\begin{proposition} The morphism $\tsG_{\sE}\to\tsG^D_{\sE}$ 
factors through $\tsG^+_{\sE}$.
\begin{proof} This is an immediate consequence of Proposition \ref{natural}
and the existence of a projection $\sF\to\sH_5$ of bundles over 
$\sE\times_{\sM_{Ell}}S$.
\end{proof}
\end{proposition}
\end{section}
\begin{section}{Subregular unstable bundles and del Pezzo surfaces}\label{subregular}
We keep the notation of Section \ref{previous}.
In particular, $r=4$ or $5$.

There is a decomposition
$${\gg}_r:=\ad\Xi_r={\ll}_r\oplus{\uu}_r\oplus{\uu}_r^*,$$
of vector bundles on $\sE$,
where $\ll_r=\ad\lambda_r$ consists of summands of degree $0$,
$\uu_r$ consists of summands of positive degree and
$\ll_r\oplus\uu_r=\ad (\lambda_r\times^{\Lambda_r}P_r)$.

We know that there
is a miniversal deformation space $\szd$ for $\Xi_5$
with a morphism $\szd\to[\hsY/W]$ that preserves $\GG_m$-orbits;
note that this is weaker than the statement that the
morphism is $\GG_m$-equivariant.

Define $\tszd=\szd\times_{\sG_{\sE}}\tsG_{\sE}$, with projection
$\nu:\tszd\to\szd$. By restricting to $\sG^{ss}_{\sE}$ 
we see that $\nu$ is generically finite
and that $\deg\nu=\#W$.
The existence of the morphisms
$\pi:\tsG_{\sE}\to \sY$ and $\szd\to[\hsY/W]$ gives a square
$$\xymatrix{
{\tszd}\ar[r]\ar[d] & {\sY}\ar[d]\\
{[\hsY/W]}\ar@{-->}[r] & {[\sY/W]}
}$$
where the lower broken arrow is a rational map,
the projection of a cone from
its vertex to its base. On the open substack $\nu^{-1}(\szdss)$
this square is commutative, and so there is a factorization
through the line bundle $\sL\to \sY$ corresponding to the cone $\hsY$:
$$\xymatrix{
&{\tszd}\ar[dl]_{\nu^\dagger}\ar[d]\ar[dr]&\\
{\szd}\ar[d]&{\sL}\ar[r]\ar[d] & {\sY}\\
{[\hsY/W]}&{\hsY}\ar[l]&
}$$
where $\deg\nu^\dagger=\#W$.

The next lemma is a slight strengthening of \ref{smoothness}.

\begin{lemma}\label{smoothness plus}
$\tszd\to\Pic^1(\sE)\times_{\sM_{Ell}}\sC_{pre}\times_{\sM_{Ell}}\sY$
is smooth.
\begin{proof} $\tszd$ classifies $B$-bundles over pre-stable curves.
The obstruction to smoothness lies in an extension of groups
each of the form $H^2(C,\sU)$ where $C$ is a pre-stable
curve and $\sU$ is a commutative unipotent group scheme over $C$.
Since $C$ is $1$-dimensional this obstruction vanishes.
\end{proof}
\end{lemma}

Now restrict to the point $[0_{\sE}]$ of $\Pic^1(\sE)$;
that is, turn off $\Pic^1(\sE)$. Let $\sZ$ and $\tsZ$
denote the resulting spaces; they fit into a commutative diagram
$$\xymatrix{
&{\tsZ}\ar[dl]_{\nu}\ar[d]^{\mu}\ar[dr]&\\
{\sZ}\ar[d]&{\sL}\ar[r]\ar[d] & {\sY}\\
{[\hsY/W]}&{\hsY}\ar[l]&
}$$
where also $\deg\nu=\#W$.
Note that, by Lemma \ref{smoothness plus}, 
$\tsZ\to \sC_{\sM_{Ell}}\sY$ is smooth and the unstable locus $\tsZ^u$ in $\tsZ$
(defined as the inverse image $\nu^{-1}(\sZ^u)$ of the unstable
locus $\sZ^u$ in $\sZ$)
equals $\mu^{-1}(0_{\sL})$, where $0_{\sL}$ is the zero-section of $\sL$.
Moreover, $\tsZ$ is smooth over the stack $\sC_{pre}$
of elliptic curves with rational tails, and $\tsZ^u$ is the inverse
image of the discriminant divisor $\Delta_{\sC_{pre}}$ in $\sC_{pre}$.
So $\mu^{-1}(0_{\sL})$ has normal crossings. However, we do not yet know that
$\mu^{-1}(0_{\sL})$ is reduced; this will be proved in Theorem \ref{7.3} below.

By construction, $\sZ^u$ is the inverse image of the vertex of the cone
$[\hsY/W]$, so that $\sZ^u$ is a surface, relative to $\sM_{Ell}$,
with two strata:
$\sZ^u=\{[\Xi_5]\}\cup(\sZ^u-\{[\Xi_5]\})$; the points of
$\sZ^u-\{[\Xi_5]\}$ correspond to regular unstable bundles, all of which
are isomorphic to $\Xi_4$.

The next result, giving a version of Steinberg's
cross-section theorem which underlies the classical
BGSS construction, is
due to Friedman and Morgan [FM], Theorem 5.1.1;
see also [Brue]. From it, they
get a new proof of Looijenga's result [L1],
that $[\hY/W]$ is an affine space.
They all assume that the ground field is $\C$
but specialization extends the result to all characteristics,
as we show.

\begin{theorem}\label{FM} 
\part[i] [FM] For any geometric point $\eta$ in $\sZ^u-\{[\Xi_5]\}$
there is a section of the morphism $\sZ\to[\hsY/W]$ through $\eta$.
\part[ii] [L1] $[\hsY/W]$ is isomorphic to an affine space 
bundle over $\sM_{Ell}$ whose fibre is $\A^{l+1}$.
\begin{proof} Friedman and Morgan
construct, for the regular
unstable bundle $\Xi_4$, a chart $V$ that is isomorphic to 
$H^1(\sE,\ad\Xi_4)$, a vector bundle over $\sM_{Ell}$. 
Their main result is that
the classifying morphism $\hat{h}:V\to [\hsY/W]$, which preserves
$\GG_m$-orbits, induces a morphism $h:\P(V)\to [\sY/W]$
defied over $\sM_{Ell}$ which is an isomorphism 
over $\sM_{Ell,\Q}=\sM_{Ell}\otimes_{\Z}\Q$.

Now the formation of geometric quotients
does not always commute with reduction modulo $p$,
so we know only that $[\sY/W]$ is normal and the geometric fibres 
of $[\sY/W]\to\sM_{Ell}$ are generically reduced
(on the locus where $W$ acts freely).
However, $h:\P(V)\to [\sY/W]$ provides a simultaneous normalization
of the fibres of the family 
$[\sY/W]\to\sM_{Ell}$ (without the requirement for any base change)
and it follows that $h$ is an isomorphism.

Moreover, the weights of the $\GG_m$-action on $V$ and on $[\hY/W]$
are the same on both sides, so that $h^*L\cong\sO_{\P(V)}(1)$
and $\hat{h}$ is an isomorphism.

By specialization, it follows that $h$ is an isomorphism in 
characteristic $p$, and then we deduce similarly that $\hat{h}$
is an isomorphism in characteristic $p$.

By the openness of versality, $\sZ$ is also a chart for $\eta$; that is,
$\sZ\to\sG_{\sE}$ is smooth at $\eta$, 
so an appropriate slice of $\sZ$ through $\eta$
gives the section required.
\end{proof}
\end{theorem}

Now fix a geometric point of $\sM_{Ell}$,
corresponding to an elliptic curve $E$
over an algebraically closed field $k$.
Fix the origin $0_Y$ of $Y=\Hom(P,E)$ and the copy of $\A^1_k$ that is the 
line in $L=\sL\vert_Y$ lying over $0_Y$. Let $\tX\to\A^1_k$
denote the restriction of $\tZ$ to this line. That is,
$\tX$ is the fibre over $0_Y$ of the smooth morphism
$\tZ\to Y$, so that $\tX$ is smooth, of dimension 
$\dim\tX=\dim\tZ-\dim Y=l+3-l=3$, and the zero fibre
$\tX_0$ of $\tX\to\A^1_k$ has normal crossings.

Now $\A^1_k$ maps isomorphically to its image $\G$
in $[\hY/W]$, since $W$ acts trivially on this line $\A^1_k$,
and so there is a commutative diagram
$$\xymatrix{
{\tX}\ar[d]_f\ar[dr]\\
{X}\ar[r]&{\A^1_k}
}$$
where $X$ is the image of $\tX$ in $Z$ and $\tX\to X$ is proper and birational.
Moreover, this diagram is $\GG_m$-equivariant,
since $f$ takes $\GG_m$-orbits to $\GG_m$-orbits and is birational,
$0\in\A^1_k$ is a fixed point and $\A^1_k-\{0\}$ is homogeneous
(but not necessarily a torsor; this will be proved later,
in Lemma \ref{6.16plus}) under the $\GG_m$-action.

\begin{lemma} $X$ has local complete intersection (LCI) singularities
and $X_0$ is a normal surface with a unique singular point, namely
the point $[\xi]$ corresponding to $\xi$.
\begin{proof} $X$ is the inverse image in $Z$ of the line $\G$ in $[\hY/W]$.
Since $[\hY/W]\cong\A^{l+1}_k$ and $Z$ is smooth, the singularities of $X$
are LCI. That $[\xi]$ is isolated follows from the fact that for the exceptional
groups there is (up to choosing a point on $\Pic^1(E)$, which amounts
to translating $\xi$ by a point on $E$) just one subregular unstable bundle
and the subregular unstable locus has
codimension $2$ in the unstable locus.
\end{proof}
\end{lemma}

Write $\sF_r=\Xi_r\times^GF$, the associated $F$-bundle over $E$. We know,
from the computation of Helmke and Slodowy, that
$\dim H^0(E,\gg_5)=l+4$.

\begin{proposition}\label{negative}\label{7.8}\label{6.4} 
Suppose that $\s$ is a section of $\sF_r$ 
that gives a negative cocharacter $[\s]$.

\part[i] If $r=4$, then $[\s]=-\a_4^\vee$.
\part[ii] If $r=5$, then $[\s]=-\a_4^\vee$ or
$-\a_5^\vee$ or $-\a_4^\vee-\a_5^\vee$.
\begin{proof}
A section $\s$ of $\sF_r$ corresponds to a $B$-bundle $\sB_\s\to E$
such that the induced bundle $\sB_\s\times^BG$ equals $\Xi_r$.
  
Suppose $\mathfrak u_\Sigma$ is an ideal in $\Lie B$,
$\Sigma$ the corresponding subset of the set of positive roots, 
$\uu_\Sigma = \mathfrak u_\Sigma \times^B \sB$ the 
bundle attached to this by $\sB$ 
and $2\rho_\Sigma = \sum_{\alpha \in \Sigma} \alpha$.
Then, writing $[\s]=-\sum_{i\in I}s_i\a_i^\vee$, $s_i\ge 0$, we have 
$$ \deg \uu_\Sigma = \left(-[\s],2\rho_\Sigma\right) = 
\left(\sum s_i\a_i^\vee,
2\rho_\Sigma\right).$$
As $\uu_\Sigma$ is a sub-bundle of $\gg$, we have
 $\dim H^0(E,\gg_r)\ge\dim H^0(E,\uu_\Sigma)\ge\deg\uu_\Sigma$.

\begin{lemma}
$\dim H^0(E,\gg_r) > \deg\uu_\Sigma$. 
\begin{proof}
Suppose otherwise. Then $\deg \uu_\Sigma =
\dim H^0(E,\uu_\Sigma)$, and so 
the only possible indecomposable summands of the vector bundle
$\uu_\Sigma$  with slope zero must have non-trivial determinant; the
remaining summands have strictly positive slope.

Next, consider the 
canonical reduction of $\gg_r$.
This is $\gg_r=\ll_r + \nn^+ + \nn^-$, where 
$\ll_r =\zeta({\ll_r}) + [\ll_r,\ll_r]$ is the decomposition
of $\ll_r$ as the direct sum of its centre and the bundle of
its derived subalgebra, 
$\zeta({\ll_r})$ is a non-zero direct sum of trivial line
bundles, $[\ll_r,\ll_r]$ is a direct sum of semi-stable bundles of degree 0,
and $\nn^+$
(resp.{} $\nn^-$) is a direct sum of vector bundles of strictly positive
(resp.{} negative) slopes. Then $\uu_\Sigma $ embeds into
$[\ll_r,\ll_r] +\nn^+$, and hence 
$\dim H^0(E,\uu_\Sigma) \leq \dim H^0(E,[\ll_r,\ll_r]+\nn^+) 
< \dim H^0(E,\gg_r)$.
\end{proof}
\end{lemma}

Now let $J$ be a subset of the simple roots, and take $\Sigma_J$ to be
the roots of the unipotent radical of the standard parabolic $P_J$ generated
by $B$ and the negative simple roots $-\alpha$  for $\alpha \not\in J$, 
and write $\rho_J$ for $\rho_{\Sigma_j}$, $\uu_J$ for $\uu_{\Sigma_J}$.

Suppose that $\a$ is simple and $\a\not\in J$. Let ${\mathfrak s}_\a$ denote
the reflection in $\a$. 
Then ${\mathfrak s}_\alpha(\Sigma_J) =\Sigma_J$ and
${\mathfrak s}_\alpha (\alpha^\vee) = - \alpha^\vee$, so $(\alpha^\vee, \rho_J) =
0$, and hence $2\rho_J = \sum_{j\in J}  m_j\varpi_j$ for some integers 
$m_j$ depending on $J$. In particular, $2\rho_{\{j\}} = m_j \varpi_j$
and $2\rho_{\{j,k\} } = a_j\varpi_j + a_k \varpi_k$ for some integers
$m_j$, $a_j$, $a_k$.
Evaluating on 
$\varpi_j^\vee$, we get 
$ m_j(\varpi_j^\vee, \varpi_j ) = (\varpi_j^\vee,2\rho_{\{j\}} ) = 
(\varpi_j^\vee,2\rho) $,
and so the numbers $m_j$ can be read off
from Bourbaki's planches; the vector $m=(m_1,...,m_l)$
is given as follows, and 
$\dim H^0(E,\gg_r) > \deg \uu_{\{j\}} = s_jm_j$ for all $j$.
\begin{eqnarray*}
E_8:\  m&=&(23,17,13,9,11,14,19,29);\\
E_7:\  m&=&(17,14,11,8,10,13,18);\\
E_6:\  m&=&(12,11,9,7,9,12);\\
D_5:\  m&=&(8,7,8,6,8).
\end{eqnarray*}

\DHrefpart{i}: $r=4$. Then
$l+2 = \dim H^0(E,\gg_4), $ 
$s_j  = 0$ if $j \neq 4$, and $s_4 \leq 1$. 
But $[\s]<0$, so
$s_4 > 0$ and then
$s_4 = 1$. This proves \DHrefpart{i} of the Proposition.

\DHrefpart{ii}: $r=5$. 
$\Xi_5$ is
induced from a semi-stable $\Lambda_5$-bundle 
[HS2], and 
$l + 4 = \dim H^0(E,\gg_5) $, 
so again $s_j \le 1$.

To finish, we need to show
that if $j \neq 4$ or $5$, then $s_j = 0$.
If $G=E_7$ or $E_8$, this follows from 
$l+4=\dim H^0(E,\gg)>\dim H^0(E,\uu_{\{j\}})=s_jm_j$
and the values given above for the vector $m$.
For $G=E_6$ or $D_5$ we can (eschewing science) calculate $\deg \uu_{\{ 5,j\}} $
case by case, as necessary. If $G=E_6$ and $s_3 = 1$, we get 
$\deg \uu_{\{5,j\}} =10 = l+4$, 
while for $G=D_5$ we get $\deg \uu_{\{5,j\}} =11,9,10$ when
$j =1,2,3$, which furnishes a contradiction
and thereby concludes the proof of the Proposition.
\end{proof}
\end{proposition}

We remark that a uniform treatment of all the combinatorics in this
paper will be found in [GS].

\begin{corollary} The morphism $p:\tsG_{\sE}\to\sG_{\sE}$ is representable
over a neighbourhood of $\Xi_5$ that contains $\sG^{ss}_{\sE}$.
\begin{proof} From Proposition \ref{7.8}, any rational component of a stable
map $C\to\sF_5$ maps birationally to its image if
it is not contracted to a point. So the map $C\to\sF_5$ has no automorphisms.
\end{proof}
\end{corollary}

The next result is the crux of our paper. It is this that shows
how exceptional groups lead naturally to weak del Pezzo surfaces.
The strategy of the proof is to show first, by considering which components
can meet, that the only triple points on the closed fibre $\tX_0$ lie on certain
copies $Q$ of $\P^1\times\P^1$ 
and then to compute intersection
numbers in as many ways as possible and use the adjunction
formula as much as possible in order to derive the full result.

Note that, although the surface $Q$ that appears might appear
to be annoying, its presence enables us to recognize
the other surface $D_1$ as a weak del Pezzo surface, because
the intersection $D_1\cap Q$ is a line on $D_1$: if $Q$, and so that
line, were missing then
our argument would only show that $D_1$ is a surface whose first 
Chern class has a sign (modulo $(-2)$-curves).
In the course of the proof we shall see that there is
a smooth birational contraction where
$Q$ is contracted onto that line.

\begin{theorem}\label{7.3}\label{sst} 

\part[i] The fibre $\tX_0$ is semi-stable.
That is, it is a reduced union of smooth surfaces and it has normal crossings. 

\part[ii] $\tX_0$ consists
of three components: $\tX_0=D_0+D_1+Q$, where $D_0$ is the strict transform
of the affine surface $X_0$, 
$Q$ is a copy of $\P^1\times\P^1$, the normal bundle
$\sN_{Q/\tX}$ is isomorphic to $\sO(-1,-1)$ and $D_1$ is a 
weak del Pezzo surface.

\part[iii] Set $D_0\cap D_1=A$. Then $A$ is isomorphic to $E$
and $A$ is an anticanonical divisor on both $D_0$ and on $D_1$.

\part[iv] $Q\cap D_0$ is a fibre $\beta$ of one ruling on $Q$,
$Q\cap D_1$ is a fibre $\gamma$ of the other ruling, $\gamma$
is the unique $(-1)$-curve on $D_1$
and $\beta$ is a $(-1)$-curve on $D_0$.

\part[v]\label{chain}
There is a number $t$ with $0\le t\le l-1$ that depends
only upon $G$ such that 
the dual graph of the
exceptional locus of the contraction $D_0\to X_0$
is a chain
$$\xymatrix{
*+[Fo]{\epsilon_{t}}\ar@{-}[r]&{}\ar@{.}[r]&{}\ar@{-}[r]&  
*+[Fo]{\epsilon_1}\ar@{-}[r]&*+[Fo]{\b}\ar@{-}[r]&*+[Fo]{A}            
}$$
of transverse curves where each $\epsilon_i$ is a $(-2)$-curve.

\part[vi] There is a birational map
$\tX-\to X^-$, where $X^-$ is smooth and
$X^-\to\A^1$ is also semi-stable, 
which is constructed by first contracting $Q$
onto $\g$ and then successively flopping the strict transforms of
$\epsilon_1,...,\epsilon_t$.

\part[vii] 
$X^-_0=D_0^-+D_1^-$ where $D_0^-$ is the minimal
resolution of $X_0$, $D_1^-$ is a weak del Pezzo surface
of degree $9-l$ and $D_1^-$ contains a configuration
$\G_0$ of $(-2)$-curves of type $G$.
\begin{remark} We do not need to know the value of $t$,
since the curves $\b,\epsilon_i$ will be flopped
to give a family of del Pezzo surfaces of the correct degree.
This issue will be addressed in Davis' thesis [D].
\end{remark}
\begin{proof}
The points $x$ of $\tX_0$ correspond to configurations
$C=\s+\phi$, where $\s$ is a section of either $\xi\times^GF$ or
$\eta\times^GF$, where $\eta$ is regular unstable.
That is, $\eta\in X_0-\{\xi\}$. 

Since $X_0$ is reduced, by Theorem \ref{FM}, it follows 
from Proposition \ref{negative} that
$$\tX_0=D_0+\sum_{i=1}^r m_iD_i+\sum_{j=1}^s n_jQ_j +\sum r_kF_k,$$
where $D_0$ is the strict transform of $X_0$ and the other components
are projective. Moreover,
points in the interior of the various components correspond
to configurations $\s+\phi$ where $\s$ is a copy of $E$ and the class of
$\phi$ is given as follows:
\begin{enumerate}
\item $D_0$: $[\phi]=\a_4^\vee$ or $\a_5^\vee$;

\item $D_i^o$ for $i\ge 1$: $[\phi]=\a_5^\vee$;

\item $F_k^o$: $[\phi]=\a_4^\vee$;

\item $Q_j^o$: $[\phi]=\a_4^\vee+\a_5^\vee$.
\end{enumerate}

\begin{lemma}\label{disjoint} For $j\ge 1$ the $D_j$
are smooth and are mutually disjoint,
as are the $Q_{j}$ and the $F_k$.
\begin{proof} If two branches of the divisor
$\sum_{j>0}D_j$ were to meet, then they would meet along a curve
$\delta$ that parametrized
sections $\tau$ of $\sF_{\xi}$ with
$[\tau]=-\a_4^\vee-\a_5^\vee$; more precisely, $\delta$
would parametrize
stable maps whose image is $\tau+\phi_4+\phi_5$ with $[\phi_i]=\a_i^\vee$.
Since $\tsG\to\sC_{pre}$ is smooth, the locus of stable maps whose image contains
a curve of class $\a_4$ is of pure codimension; however, it contains
$\delta$ and the surface $D_0$. So $\sum_{j>0}D_j$ is smooth.

The argument for $\sum Q_j$ and $\sum F_k$ is similar.
\end{proof}
\end{lemma}

\begin{lemma}\label{6.9}
$Q_j$ is isomorphic to $\P^1\times\P^1$.
\begin{proof}
Every point of $Q_j$ corresponds to a stable map
$f:C=\s\cup \phi\to\sF_\xi$ where $\s\cap \phi$ consists of a single point,
say $x$, and $[\phi]=\a_4^\vee+\a_5^\vee$.
Therefore $Q_j$ is isomorphic to 
the variety $V$
of $(1,1)$ curves on the flag variety $H=SL_3/B$ of type $A_2$
that pass through a fixed
point $v_0$ of $H$
($v_0$ is also here
a fixed point on $\s$).

So it is enough to prove that $V\cong\P^1\times\P^1$.

Recall that $H$ is the incidence variety
$\{(x,l)\vert x\in l\}$, where $x$ is a point in $\P^2$
and $l$ is a line. The projections
$p:H\to\P^2$, $q:H\to (\P^2)^\vee$
map $V$ to the set of points on a fixed line
and to the set of lines through a fixed point, respectively.
This gives a morphism $a:V\to\P^1\times\P^1$.

Suppose that $\g,\delta$ are $(1,1)$ curves on $H$
with $p(\g)=p(\delta)=m$ and $q(\g)=q(\delta)=n$.
So $m,n$ are lines. Then $\g,\delta$ are curves
in $p^{-1}(m)$ with equal image under $q$.
Now $p^{-1}(m)$ is a copy of the Hirzebruch surface $\F_1$
and $q:\F_1\to(\P^2)^\vee$ exhibits $\F_1$
as the blow-up of $(\P^2)^\vee$ at a point.
Then $\g=\delta$, so that $a$ is an isomorphism,
as required.
\end{proof}
\end{lemma}

\begin{lemma}
$Q_j\cap (\tX_0-Q_j)=\b\cup\g$ where $\b,\g$
are opposite rulings on $\P^1\times\P^1$.
\begin{proof}
$Q_j\cap (\tX_0-Q_j)$ parametrizes
stable maps which are \emph{either} of the form
$f:C=\s\cup\phi_4\cup\phi_5\to \sF_\xi$
\emph{or} of the form
$f:C=\s\cup\psi\cup\phi_4\cup\phi_5\to \sF_\xi$
where in each case $[f(\phi_i]=\a_i^\vee$
and $f(\psi)$ is a point.
Also $\s\cap(\phi_4\cup\phi_5)$ is a single point,
resp., $\s\cap(\psi\cup\phi_4\cup\phi_5)$ is a single point.
There are three possibilities.
\begin{enumerate}
\item $\s$ and $\phi_4$ remain fixed and $\phi_5$ moves along 
$\phi_4$. This gives $\b$.

\item $\s$ and $\phi_5$ remain fixed while $\phi_4$ moves along
$\phi_5$. This gives $\g$.

\item All remain fixed. This gives $\b\cap\g$.
\end{enumerate}
\end{proof}
\end{lemma}

\begin{lemma} Each triple point of $(\tX_0)_{red}$ lies in
some $Q_j$.
\begin{proof} A triple point is a stable map $f:C\to\sF_xi$
where $C$ has at least three double points.
The only possibility is $C=\s\cup\psi\cup\phi_4\cup\phi_5,$
which corresponds to a point $\b\cap\g$ as in the 
proof of the previous lemma.
\end{proof}
\end{lemma}

We now quote what seems to be folklore.

\begin{theorem}\label{folklore}
$\sM_{Ell}$ is simply connected.
\begin{proof} This can be found at

\noindent
https://mathoverflow.net/questions/105047/fundamental-group-of-the-moduli-stack-of-elliptic-curves
\end{proof}
\end{theorem}

\begin{lemma} Every double curve of $\tX_0$ that does not lie
in $\sum Q_j$ meets at least one $Q_j$ and dominates $E$.
\begin{proof} We can work with the universal curve $\sE$.
Since $\sM_{Ell}$ is simply connected, every geometric
component of $\tX_0$ and every double curve
is defined over $\sM_{Ell}$.

Suppose that $A$ is a double curve meeting no $Q_j$.
Then $A$ parametrizes stable maps $f:\s\cup\phi_4+\phi_5\to\sF_\xi$
where $\s\cap\phi_4$ is never equal to $\s\cap\phi_5$.
This defines a subtraction morphism $t:A\to E$ by
$t(f)=(\s\cap\phi_4)-(\s\cap\phi_5)$ which is never zero,
and so is constant and defines a non-zero point on $E$.
But on the universal curve the only
rational point is zero.
\end{proof}
\end{lemma}

\begin{lemma}\label{6.14}
$\tX_0$ is reduced and $K_{\tX}\sim\sum Q_j$.
\begin{proof}
We argue by induction.

Start by considering a neighbourhood in $\tX$
of some $Q_j=Q$. It meets two components $L, M$ of $\tX_0$.
Assume, as our induction hypothesis, that $L$ has multiplicity
$1$ in $\tX_0$. Since $\tX_0\sim 0$ as divisors on $\tX$,
we can write
$$\tX_0=L+mM+qQ,\ K_{\tX}\sim pM+rQ,$$
where $m,q,p,r\in\Z$ and $m,q>0$.
Say $Q\cap M=\g$, $Q\cap L=\b$.
Note that $L.\b=(\b_Q)^2=0$ and $M.\g=(\g_Q)^2=0$.

So $\tX_0.\b=0$ gives $m=qa$
and $\tX_0.\g=0$ gives
$qb=1$, so that $q=b=1$ and $m=a$.
 
The adjunction formula for the curve $\g$ in the $3$-fold $\tX$
is $2g(\g)-2=K_{\tX}.\g+\deg\sN_{\g/\tX}$, which 
gives $-2=r.(-1)-1$, so $r=1$ and we have
$$\tX_0=L+mM+Q,\ K_{\tX}\sim pM+Q.$$

The adjunction formula applied to $\b$ gives
$p=2(m-1)$, so that $K_{\tX}\sim 2(m-1)M+Q$.

Suppose next that $A$ is a double curve
that lies in $L$ and meets $Q$. So $A$ is a connected component of
$L\cap M$, meeting
$Q_1,...,Q_s$, say.
Pick a sufficiently small neighbourhood $U$ of $A$ that contains $Q_1,...,Q_s$.

The adjunction formula for $L$ gives
$K_L\sim (m-2)A$ (in the neighbourhood $U\cap L$
of $A$ in $L$).
There are $(-1)$-curves $\b_1,...,\b_s$ in this neighbourhood;
they satisfy $\b_i.A=1$ and $\b_i.K_L=-1$,
and so $m=1$ and then $K_{\tX}\sim Q$.
Since $\tX_0$ is connected in codimension one,
this completes the proof of Lemma \ref{6.14}.
\end{proof}
\end{lemma}

Now drop the distinction between the $D_i$ and the $F_k$.
The adjunction formula shows that
$-K_{D_i}\sim\left(\sum_{j\ne i}D_j\right)\vert_{D_i}$.

\begin{lemma} If $L$ is a smooth projective surface
and $B\in\vert-K+L\vert$ is smooth, but maybe disconnected,
then there are three possibilities.

\begin{enumerate}
\item $B$ is connected, $g(B)=1$ and $L$ is rational.

\item $B$ is connected, $g(B)=1$, $L$ is elliptic ruled
and $B$ is a bisection of the ruling.

\item $B$ has just two connected components $B_1$ and $B_2$,
$g(B_i)=1$, $L$ is elliptic ruled
and each $B_i$ is a section of the ruling.
\end{enumerate}
\begin{proof} This is easy and well known.
\end{proof}
\end{lemma}

Since $X_0$ has a $\GG_m$-action, $D_0$ contains at most one curve
of genus $\ge 1$. 
Since $D_0$ meets $\tX_0-D_0$, the divisors
$D_0,...,D_r$ form a chain
of surfaces,
where, for $r\ge i\ge 1$, each $D_i\cap D_{i-1}$ is a curve $A_i$
of genus $1$. Moreover, $A_1$ is an anticanonical divisor on
$D_0$, 
$A_{r}$ is an anticanonical divisor on $D_r$ and,
for $1\le i\le r-1$,
$A_i+A_{i+1}$ is an anticanonical divisor on $D_i$.
So $D_1,...,D_{r-1}$ are elliptic ruled and
$D_r$ is rational or elliptic ruled.

\begin{lemma} Each $A_i$ meets exactly one $Q_j$ and is isomorphic to $E$.
\begin{proof} 
If $A_i$ meets $s$ of the $Q_j$ then the subtraction morphism
$t:A_i\to E$ has a fibre $t^{-1}(0)$ consisting of $s$ points.
So $t$ is {\'e}tale of degree $s$ in characteristic zero.
Then $t$ is {\'e}tale of degree $s$ universally. Then the kernel
of the dual
morphism $t^\vee:E\to \Pic^0(A_i)$
is a finite {\'e}tale subgroup scheme of
$E$ of order $s$. But on the universal curve
the only such subgroup scheme is trivial,
so $s=1$.
\end{proof}
\end{lemma}

So $\tX_0=\sum_0^rD_i+\sum_1^rQ_i$
where the $D_i$ form a chain of surfaces
that intersect in copies $A_1,...,A_r$ of $E$
and $Q_j$ intersects $D_{j-1}$ in the curve $\b_j$
and $D_j$ in the curve $\g_j$. Both $\b_{j+1}$
and $\g_j$ are $(-1)$-curves on $D_j$.
The dual complex of $\tX_0$ is shown in Figure \ref{fig:dual complex}.
\begin{figure*}[htbp]
\centerline
{\includegraphics[height=0.5\textheight, width=0.5\textwidth, angle=270]{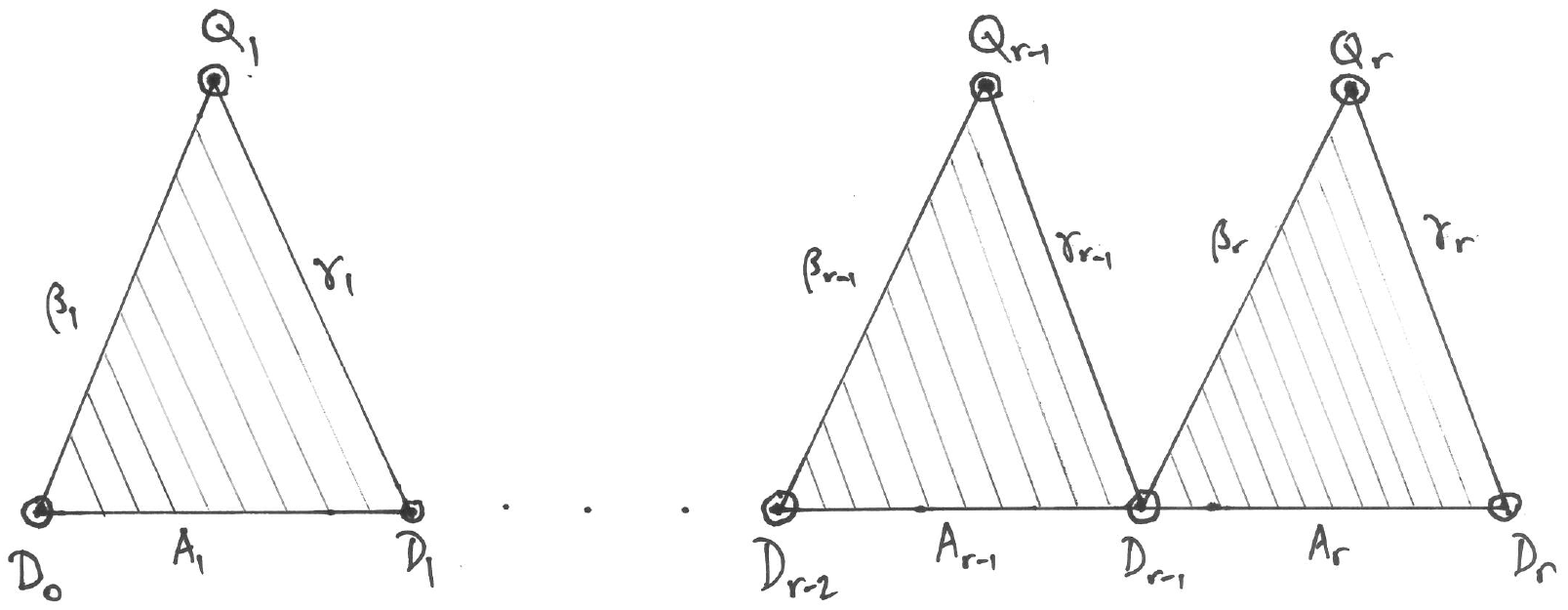}}
\caption{The dual complex of $\tX_0$}
\label{fig:dual complex}
\end{figure*}

\begin{lemma}
$D_r$ is geometrically rational.
\begin{proof} If not, then $A_r$ is a bisection
of the ruling on $D_r$, so that there is a universally defined
morphism $E\to E$ of degree $2$, which is impossible.
\end{proof}
\end{lemma}

\begin{corollary}\label{HS redux} $X_0$ has a simply elliptic singularity.
Say its degree is $d$; then
$d\le 4$. Its multiplicity is $\max\{2,d\}$ and its embedding
dimension is $\max\{3,d\}$.
\begin{proof} 
We know that $K_{D_0}\sim-A_1$;
it follows at once that $X_0$ is simply elliptic.
That is, the exceptional locus in
its minimal resolution is a single elliptic curve.
Since $X_0$ 
has a $\GG_m$-action
it is isomorphic to a line bundle over $A_1$,
which is isomorphic to $E$.

There is a Cartesian diagram
$$\xymatrix{
{X_0}\ar@{^(->}[r]\ar[d]&{X}\ar@{^(->}[r]\ar[d]& {Z}\ar[d]\\
{\{0\}}\ar@{^(->}[r]&{\A^1}\ar@{^(->}[r] & {\A^{l+1}}
}$$
whose horizontal arrows are closed embeddings; it
exhibits $X_0$ as an LCI singularity.
From the well known structure of these 
homogeneous co-ordinate rings
it follows that $d\le 4$
and that the multiplicity and embedding
dimension of $X_0$ are as described. 
\end{proof}
\end{corollary}

We now let $d$ denote
the degree of the simply elliptic
singularity $X_0$.
Recall that $l$ is the rank of the group $G$.

\begin{proposition}\label{blank} $d=9-l$.
\begin{proof} The point is to show that $X_0$ has, first, the correct 
embedding dimension and, then, the correct multiplicity.

Recall that $Z\cong\A^{l+3}$ and that $X\to\A^1$ is obtained
from a morphism $\pi:Z\to[\hY/W]$ that preserves $\GG_m$-orbits by restricting
to the line in $[\hY/W]$ that goes through the origin of $Y$.
So $X_0$ is embedded in $Z$ as the fibre over the vertex of
$[\hY/W]$. By Looijenga's result
(Theorem \ref{FM} above), $[\hY/W]$ is isomorphic to $\A^{l+1}$
and the weights of the $\GG_m$-action
on $[\hY/W]$ are the coefficients
in the biggest root in the affine Dynkin diagram.
These are given by Table $1$,
in which each exponent denotes the multiplicity of the relevant weight.
\begin{table}[ht]
\centering 
\caption{Coefficients in the biggest root in the affine diagram}
$$\begin{array}{c c c c c c c} 
{E_8:} & 1 & {2^2} & {3^2} & {4^2} & 5 & 6\\
{E_7:} & {1^2} & {2^3} & {3^2} & 4 & {} & {}\\
{E_6:} & {1^3} & {2^3} & 3 & {} & {} & {}\\
{D_5:} & {1^4} & {2^2} & {} & {} & {} & {}\\
\end{array}$$
\end{table}

We also need to know the weights of the $\GG_m$-action coming
from the action of the centre $\zeta(\Lambda_5)$ of $\Lambda_5$.
Note that $\zeta(\Lambda_5)\cong\GG_m$.
\begin{lemma}\label{weights}\label{7.12}
The weights of the $\zeta(\Lambda_5)$-action on $Z$ are given in Table $2$.
\begin{table}[ht]
\centering
\caption{Weight multiplicities for the action of $\zeta(\Lambda_5)$ on $Z$}
$$
\begin{array}{c c c c c c} 
{E_8:} & {1^2} & {2^3} & {3^3} & {4^2} & 5\\
{E_7:} & {1^4} & {2^4} & {3^2} & {} & {}\\
{E_6:} & {1^6} & {2^3} & {} & {} & {}\\
{D_5:} & {1^8} & {} & {} & {} & {}\\
\end{array}
$$
\end{table}
\begin{proof}
As before, we have $\gg_5=\ll_5\oplus\nn^+\oplus\nn^-$.
This represents 
the fact that the $G/P_5$-bundle associated to $\Xi_5$
has a section $\s$ with $[\s]=-\a_5^\vee$
(see Remark \ref{explain}).

Consider the grading of $\nn^+ = \oplus_{0<i\leq r} \nn^+_i$
defined by the action of $\GG_m$ via the cocharacter
$\varpi_5^\vee$. Then
the fibre $\mathfrak n^+_i$ of $\nn^+_i$  is the direct
sum $\mathfrak g_{\a}$ of root spaces 
where $\a$ runs over the set $I_5(i)$ of positive
roots such that $\a_5$ has coefficient exactly $i$ in $\a$, and  $r$ is 
the multiplicity of $\a_5$ in the biggest root.

Note that $\mathfrak n^+_i$ is a representation of 
$\Lambda_5$, so that $\sum_{\a\in I_5(i)}\a$ is trivial
on the derived subgroup of $\Lambda_5$. Therefore
$\sum_{\a\in I_5(i)}\a= n_i\varpi_5$, where $n_i$ depends on $i$ and $G$.

Put $d_i=\deg\nn^+_i$,
so that
$d_i=(\sum_{\a\in I_5(i)}\a).\a_5^\vee=n_i$.
But $n_i(\varpi_5.\varpi_5^\vee)=\sum_{\a\in I_5(i)}(\a.\varpi_5^\vee)=i\#I_5(i)$,
and then an inspection of Bourbaki's planches
shows that 
$r$ and the $d_i$ are given as follows:
$$
\begin{array}{ccccc}
{E_8:}&{r=5,} & {(d_1,...,d_5)}&{=}&{(2,3,3,2,1)}\\
{E_7:}&{r=3,} & {(d_1,d_2,d_3)}&{=}&{(4,4,2)}\\
{E_6:}&{r=2,} & {(d_1,d_2)}&{=}&{(6,3)}\\
{D_5:}&{r=1,} & {d_1}&{=}&{8.}
\end{array}
$$
From the description above of $\gg$,
we know that $H^1(E,\nn^+) = \oplus_i H^1(E, \nn^+_i) = 0$,
so by  Riemann--Roch $\dim H^0(E,\nn^+_i) = \chi(E,\nn^+_i) = d_i$,
and as $\zeta(\Lambda_5) = \GG_m$ acts with weight $i$ on the bundle $\nn^+_i$
and so on $H^0(E,\nn^+_i)$, the lemma is proved.
\end{proof}
\end{lemma}

\begin{remark}\label{ian's remark}
Note that $H^0(E,\gg_5) = H^0(E,\ll_5) \oplus H^0(E,\nn^+)$ and 
$\dim H^0(E,\ll_5) = 1$, since $\lambda_5$ is
determined by a $GL_p \times GL_q \times GL_s$ bundle
$(\eta_p,\eta_q,\eta_s)$ with $\det \eta_p = \det \eta_q = \det \eta_s$
of degree $1$ and $(p,q,s) = (1,5,l-4)$.
This provides another proof that $\dim H^0(E,\gg_5)=l+4$.
\end{remark}

\begin{lemma} If $f:A\to B$ is an injective $k$-algebra homomorphism 
of positively graded domains 
such that the corresponding morphism $F:\Sp B\to\Sp A$ takes $\GG_m$-orbits to
$\GG_m$-orbits, then there is an integer $e=e(f)>0$ such that
$f(A_n)\subset B_{en}$ for all $n$. That is, $f$
multiplies degrees by $e$.
\begin{proof} That $F$ preserves $\GG_m$-orbits means that for every $n$,
there exists $n'$ such that $f(A_n)\subset B_{n'}$. Suppose that $x\in A_m$
and $y\in A_n$; then
$$n\deg f(x)=\deg f(x^n)=\deg f(y^m)=m\deg f(y),$$
so that $\deg f(x)/\deg f(y) = m/n = \deg x/\deg y$,
which proves the lemma.
\end{proof}
\end{lemma}

The next result is crucial. To prove it we 
must use the fact that has now been proved, that
$X_0$ is a simply elliptic LCI singularity
and so has embedding dimension at most $4$.

\begin{proposition}\label{6.16} For each group $G=E_8,...,D_5$ the 
integer $e$ is equal to $1$.
\begin{proof} We check the cases separately. For $E_8$
the tables above show that $A,B$ are weighted
polynomial rings, as follows: $A=k[1,2^2,3^2,4^2,5,6]$ and
$B=k[1^2,2^3,3^3,4^2,5]$. Then the fact that the fibre $X_0$ over the origin 
of $\Sp B\to\Sp A$ has embedding dimension at most $4$
at its unique singular point
forces $e=1$ in this case, and $e=1$ in all the other cases
for exactly similar reasons.
\end{proof}
\end{proposition}

\begin{remark}
It is possible instead to
give a case-free proof of this proposition via the interpretation
of the line bundle $\sL$ on $\sY$ as a determinant line bundle;
we plan to return to this approach [GS].
\end{remark}

Therefore $X_0$ is embedded $\GG_m$-equivariantly
as a complete intersection in an affine space
as follows:
$$
\begin{array}{c c c c c c} 
{E_8:} & {X_0} & = & {(6)} & \inj & {\A^3(1,2,3)}\\
{E_7:} & {X_0} & = & {(4)} & \inj & {\A^3(1^2,2)}\\
{E_6:} & {X_0} & = & {(3)} & \inj & {\A^3(1^3)}\\
{D_5:} & {X_0} & = & {(2,2)} & \inj & {\A^4(1^4).}
\end{array}
$$
It follows that $X_0$ has the correct multiplicity, and
the proof of Proposition \ref{blank} is complete.
\end{proof}
\end{proposition}

\begin{lemma}\label{standard}\label{6.16plus}
The morphisms $\tX\to\A^1_k$ 
and $X\to\A^1_k$
cover the standard action 
of $\GG_m$ on $\A^1_k$.
\begin{proof}
It is enough to prove the lemma for $X\to\A^1_k$.
However, this is nothing but Proposition \ref{6.16}.
\end{proof}
\end{lemma}

Next, we recover the result whose proof was sketched by
Helmke and Slodowy [HS3] on
deformations of simply elliptic singularities.
In fact
the results of Hirokado [Hi], Corollary $4.3$
and his calculation of $T^1$ 
(the tangent space to a miniversal deformation space)
in the proof of [Hi] Theorem $4.4$
permit us to prove this in almost all characteristics.

{\bf{From now on we make the following assumption:
$\ch k=0$ if $G=D_5$ and $\ch k\ne 9-l$ if
$G=E_l$ with $l=6$ or $7$.}}
Note that when $G=E_8$ this excludes nothing.

Hirokado shows that then
simply elliptic singularities of degree at most $3$
are classified exactly as in characteristic zero
(in particular, they are quasi--homogeneous)
and that the structure of their miniversal deformation
spaces is uniform across characteristics.
In particular a miniversal deformation space
$T$ of such a singularity has a $\GG_m$-action
and there is a hypersurface $S$ in $T$ 
where the $\GG_m$-action has positive weights.
The complementary line is where the singularity
remains simply elliptic but the elliptic curve varies.

\begin{theorem}\label{HS}
The morphism $Z\to[\hY/W]$ is the positive weight part
of a $\GG_m$-equivariant miniversal deformation of $X_0$.
\begin{proof} Suppose that $V\to S$ is the positive weight part
of a $\GG_m$-equivariant miniversal deformation of $X_0$.
Certainly $Z\to[\hY/W]$ is $\GG_m$-equivariant and
of positive weight, so there is a $\GG_m$-equivariant 
morphism $a:[\hY/W]\to S$ and a
Cartesian
square
$$\xymatrix{
{Z}\ar[r]\ar[d]_f & {V}\ar[d]\\
{[\hY/W]}\ar[r]^a & {S.}
}$$
The derivative of the morphism $a:[\hY/W]\to S$ at the vertex of
$[\hY/W]$ is determined uniquely, and then $a$
is determined uniquely by its $\GG_m$-equivariance.
Inspection shows that the weights of the $\GG_m$ action on
the affine spaces $[\hY/W]$ and $S$ are the same, so that
$a$ is an isomorphism.
\end{proof}
\end{theorem}

In particular, the family $X\to\A^1$ is induced by
a unique $\GG_m$-equivariant morphism $\phi:\A^1\to S$.
Recall that the base $\A^1$ of $X\to\A^1$ is the generator in
the cone $[\hY/W]$ that passes through the origin $0_Y$ of $Y$.

\begin{proposition} \label{the corollary}
\part[i] Under the morphism $a$, the line
$\A^1=\Sp k[t]$ in $[\hY/W]$ is mapped to a line $\ell=\Sp k[x]$ in $S$ 
where for $x\ne 0$ the surface $V_x$ has a du Val singularity of
type $G$.
\part[ii] \label{configuration}
For $t\ne 0$, $\tX_t$ 
is the minimal resolution
of $X_t$; it contains a $(-2)$-configuration $\G_t$ of type $G$.
\begin{proof}
Suppose first that $k=\C$. Then, by Theorem (7.7) of [L2],
the monodromy around $V_s$ is the full Weyl group.
That is, take a small neighbourhood $U$ of $s$ in $S$
and consider the complement $U-\Delta$ of the discriminant;
then $\pi_1(U-\Delta)$ acts on the cohomology of the generic
fibre as the full Weyl group. Since $V_s$ has at most du Val
singularities, \DHrefpart{i} is proved in this case.
Specialization then proves it for all fields.

\DHrefpart{ii} is a corollary of \DHrefpart{i} and Theorem \ref{Springer}.
\end{proof}
\end{proposition}

\begin{remark}\label{it would be nice}
It would be desirable to have a conceptual proof of this 
proposition, maybe as follows.
The copy of $\A^1=\Sp k[t]$ that appears is the line through the origin
$0_\sY$ of $\sY$ so that the bundles parametrized by $\sZ\vert_{t\ne 0}$ 
are unipotent (and semi-stable).

On the other hand, general
considerations of deformation theory show that the automorphism
group of each of these bundles has dimension strictly
less than that of $\Xi_5$, which is $l+4$. Therefore they are subregular
or regular; if we could prove that subregular bundles do arise
over points of $\A^1-\{0\}$
(that is, if we could prove directly that $\Xi_5$ deforms
to a subregular unipotent semi-stable bundle)
then Theorem \ref{Springer} could be applied, assuming
that $E$ is uniformizable; the result would follow
for all $E$ by specialization.
As it is, we only know how to deduce the existence of such a subregular
deformation from Proposition \ref{the corollary}.
\end{remark}

\begin{proposition} $r=1$.
\begin{proof} Since $K_{\tX}\sim\sum Q_i$, 
the divisors $D_1,...,D_r$ appear with
discrepancy zero in the resolution
$\tX$ of the $3$-fold $X$.
They also
dominate the closed point $0$ that represents $\xi$.
So it is enough to exhibit a partial resolution
$X'\to X$ which has c-dV singularities
and in which there is only one divisor
of discrepancy zero that dominates the origin.

For this, note that
from Lemma 6.18 and the proof of Proposition 6.17
it follows that the Cartesian diagram
$$\xymatrix{
{X}\ar@{^(->}[r]\ar[d]&{Z}\ar[d]\\
{\A^1_k}\ar@{^(->}[r] &{\A^{l+1}_k}
}$$
exhibits $X$ as a homogeneous complete intersection
of codimension $l$ in some affine space.  

For the groups $E_8,E_7,E_6$ and $D_5$, $X$ is,
respectively, a homogeneous sextic
in $\Sp k[1^2,2,3]=\A^4_k$, a quartic in $\Sp k[1^2,2]$,
a cubic in $\Sp k[1^4]$ and a complete intersection of
two quadrics in $\Sp k[1^5]$.
We know that for $t\ne 0$ the fibre $X_t$
has du Val singularities, so that 
in each case the subvariety 
$\Proj\sO_X$ in the appropriate weighted
projective space
has only du Val singularities
and is a del Pezzo surface of degree $9-l$.
Therefore the weighted blow-up $X'=\Proj \sO_X[t]\to X$ 
centred at $0$
has a unique exceptional divisor, namely $\Proj\sO_X$.
Since this divisor is Cartier and has du Val singularities,
it follows that $X'$ has only c-dV
singularities and is 
is a partial resolution of the kind that we sought.
\end{proof}
\end{proposition}

So $\tX_0=D_0+D_1+Q$ where $D_1$ is rational and $A=D_0\cap D_1$
is an effective anticanonical divisor on each of $D_0,D_1$.
Say $\b=Q\cap D_0$, $\g=Q\cap D_1$. Then each of $\b,\g$
is a $(-1)$-curve.

We define a \emph{$(-1,-2)$-chain of length $s$}
on $D_0$ to be an $s$-tuple $(\delta_1,\ldots,\delta_s)$
of smooth rational curves on $D_0$
where $\delta_1^2=-1$, $\delta_i^2=-2$ for $i\ge 2$,
$\delta_1.A=1$, $\delta_i.A=0$ for $i\ge 2$
and the configuration $A,\delta_1,\ldots,\delta_s$
is a chain.

\begin{lemma}\label{-1,-2 chain}
The only complete
curves on $D_0$ are $A$
and those curves that occur in a $(-1,-2)$-chains.
Exactly one of these chains
contains $\b$.
\begin{proof} This follows from
the adjunction formula, the fact that $A\sim-K_{D_0}$
and the fact that any configuration
of complete curves on $D_0$ is negative definite.
\end{proof}
\end{lemma}

It follows that 
there is a birational contraction $\tX\to \tX^+$
where
$Q$ is contracted onto the curve $\g$ in $D_1$
and where the images in $\tX^+$ of the curves 
in the $(-1,-2)$-chains in $D_0$
can be flopped successively onto $D_1$,
say via $\tX^+-\to\tX^-$, to give $\tX^-\to\A^1$
where $\tX^-_0=D_0^-+D_1^-$ is semi-stable, $K_{\tX^-}\sim 0$,
$D^-_0$ is the strict transform of $D_0$ and is the
minimal resolution of $X_0$,
$D_1^-$ is the strict transform of $D_1$ and is
a weak del Pezzo surface on which $A$
is an anti-canonical curve and the flopped curves on $D_1^-$
form a configuration of $(-1,-2)$-chains
which is isomorphic to the configuration of
$(-1,-2)$-chains on $D_0$.
 
Moreover, these birational birational transformations
are $\GG_m$-equivariant, so that $\tX^-$ has a $\GG_m$-action
that covers the standard action on $\A^1_k$ and is free
on $D^-_0-A^-$, where $A^-$ is the strict transform of $A$.
The action covers the standard action on $\A^1_k$ and so is free
on $\tX^--D_1^-$.

Since $(A^-_{D_1^-})^2=-(A^-_{D_0^-})^2$,
it follows that $D_1^-$ is a weak del Pezzo surface
of degree $d=9-l$.

It is easy to see that if $\GG_m$ acts on a
del Pezzo surface $S$, weak or not, and preserves
a smooth member of $\vert-K_S\vert$, then the action on $S$ is trivial.
So $\GG_m$ acts freely on $\tX^--D^-_1$ and trivially on the Cartier
divisor $D^-_1$.

The next lemma and its proof are taken almost without change from a paper by
Bass and Haboush [BH].

\begin{lemma}\label{6.10}
Assume that $S$ is a normal $k$-variety and that $D$ is a complete
relatively LCI closed subscheme of $S$ that contains every complete subscheme
of $S$. Assume also that there is a $\GG_m$-action on $S$ that is trivial on $D$
and free on $S-D$ and that every closed $\GG_m$-invariant subscheme of $S$ meets $D$.
Then there is a $\GG_m$-equivariant isomorphism $S\to N_{D/S}$
that identifies $D$ with the zero section of $N_{D/S}$.
\begin{proof} The aim is to reduce this to a situation where the arguments
of [BH] can be applied.

We need to show that there is a $\GG_m$-linearized locally free sheaf
$\sA$ on $D$ such that $S\cong\Sp\Symm^*\sA$.

Choose an open cover $S=\cup_{i\in I} S_i$ by affine $\GG_m$-invariant
open subschemes $S^{(i)}=\Sp R^{(i)}$ such that each $D^{(i)}=D\cap S^{(i)}$
is, if non-empty, defined by a regular sequence of length $r=\codim(D,S)$
in $R^{(i)}$. Say $J=\{j\in I\vert D^{(j)}\ne\emptyset\}$. Then
$S-\cup_{j\in J} S^{(j)}$ is a closed $\GG_m$-invariant subscheme of $S$
that is disjoint from $D$, so is empty. So we can suppose that each $D^{(i)}$
is non-empty.

The existence of the $\GG_m$-action is equivalent to a $\Z$-grading
$R^{(i)}=\oplus_{n\in\Z} R^{(i)}_n$. Since every closed $\GG_m$-invariant
subscheme of $S$ meets $D$, the modules $R^{(i)}_n$ vanish for $n<0$.
and the ideal of $D^{(i)}$ in $S^{(i)}$
is $R^{(i)}_+=\oplus_{n> 0} R^{(i)}_n$. So $D^{(i)}=\Sp R^{(i)}_0$,
where we regard $R^{(i)}_0$ both as a subring of $R^{(i)}$ and as the quotient
$R^{(i)}/R^{(i)}_+$.

That is, there is a unique $\GG_m$-equivariant retraction $r_i:S^{(i)}\to D^{(i)}$.
Since they are unique, these retractions glue to a $\GG_m$-equivariant retraction
$r:S\to D$. Note that, as a subscheme of $S$, $D$ contains all the closed
$\GG_m$-orbits in $S$.

Now we follow [BH], pp. 474 {\emph{et seq.}} We have $S=\Sp\sR$,
$\sR=\oplus_{n\ge 0}\sR_n$, a graded sheaf of $\sO_D$-algebras.
Set $\sI=\oplus_{n>0}\sR_n$, the ideal sheaf of the subscheme $D$ of $S$,
$\sN^\vee=\sI/\sI^2$ and $Y=\Sp\Symm^*\sN^\vee$, the normal bundle.
The tautological $\GG_m$-equivariant morphism
$u:S\to Y$ fits into a $\GG_{m,D}$-equivariant commutative diagram
$$\xymatrix{                                                                              
{S}\ar[r]^{u}\ar[d]_{r}&{Y}\ar[dl]^q\\                                                    
{D.}                                                                                      
}$$
Then the argument of [BH], p. 474, applies directly
to show that $u$ is {\'e}tale.

Let $D_0\inj Y$ be the zero section of $q$.

We now claim that $u$ is finite. To see this, take $\tS$ to be
the normalization of $Y$ in the function field of $S$. Then there
is a $\GG_m$-equivariant open embedding $S\inj \tS$
and a finite dominant morphism $\nu:\tS\to Y$ extending $u$.
Say $V=\tS-S$; this is closed and $\GG_m$-invariant in $Y$.
Since $D\subset S$ and $D$ is complete, $D$ is disjoint from $V$.
Assume that $V\ne\emptyset$; then, over any affine chart of $Y$,
there is a $\GG_m$-invariant function $f$ on $\tS$ such that
$f\vert_D=0$ and $f\vert_V=1$. However, $\sO_{\tS}^{\GG_m}=\sO_D$,
so that $V=\emptyset$ and the claim is established.

Then $u^{-1}(D_0)=D\coprod D'$ with $D'$ finite and {\'e}tale over $D_0$.
The argument just given, to show that $V=\emptyset$, can now be applied to
show that $D'=\emptyset$. So $u$ has degree $1$ over $D$,
and so is of degree $1$ everywhere. This proves Lemma \ref{6.10}.
\end{proof}
\end{lemma}

\begin{corollary}\label{normal bundle}
\part[i] There are $\GG_m$-equivariant isomorphisms
$\tX^-\to N_{D^-_1/\tX^-}$ and $N_{D^-_1/\tX^-}\to \omega_{D^-_1}$.

\part[ii] $X$ is the cone over the anti-canonical model of $D^-_1$.
\begin{proof} The first isomorphism of (1) is a special case of the lemma
and the second is a consequence of the adjunction formula. (2) is an immediate
consequence.
\end{proof}
\end{corollary}

\begin{proposition}\label{rigid dP}\label{6.30} Suppose that $D$ is a weak
del Pezzo surface of degree $9-l$ and contains a $(-2)$-configuration
$\Delta$ of type $G$. Then

\noindent (1) $D$ contains a unique line $L$;

\noindent (2) if $\d_1,...,\d_l$ are the irreducible components
of $\Delta$, numbered as in Bourbaki's planches, then $L$
meets $\Delta$ in $\d_l$ and in no other component;

\noindent (3) $\{\d_1,...,\d_l,L\}$ is a $\Z$-basis of $\NS(D)$;

\noindent (4) given a smooth member $A$ of $\vert -K_D\vert$,
$D$ is obtained by embedding $A$ in $\P^2$ via $\vert 3[0_A]\vert$
and then making $l$ successive blow-ups,
with the centre of each
blow-up being the origin $0_A$ on the strict transform of $A$.
\begin{proof} This is well known, but we include a proof for lack of
a convenient reference.

Existence of a line: $D$ is a specialization of a smooth del Pezzo $D_t$
where $-K_{D_t}$ is ample, which certainly contains lines. So $D$ does also.

Uniqueness of the line: since $\rk\NS(D)=\rk\D +1$, $\NS(D)_{\Q}$ is spanned by $L$
and $\D$ for any line $L$. Moreover, $L.\d\ge 1$ for any positive root
$\d$ with $\Supp\d=\D$. Put $H=-K_D$ and consider cases separately.

\noindent $d=1$. Then every simple root has multiplicity $\ge 2$
in the biggest root $\d_{max}$, so $L.\d_{max}\ge 2$. Then
$H^2=1=H.(L+\d_{max})$ and $(L+\d_{max})^2\ge 1$, so that, by the
index theorem, $L+\d_{max}\sim H\sim M+\d_{max}$ for any lines
$L,M$. So $L\sim M$ and then $L=M$.

\noindent $d\ge 2$. Then $(L+\d)^2\le 0$ for any line $L$ and any effective
root $\d$, by the index theorem, since $H.(L_\d)=1$ and $H^2\ge 2$.
Suppose that $\d,\e$ are effective roots with $\Supp\d=\Supp\e=\D$
and that $L,M$ are distinct lines; then $L.\d,\ M.\d\ge 1$ and
$(L+M+\d)^2\ge 2$. But $H.(L+M+\d)=2\le H^2$, so that (index)
$L+M+\d\sim H\sim L+M+\e$ and $\d\sim\e$, which is absurd,
and uniqueness is established.

Moreover, the line $L$ meets $\D$ in a simple root $\d$ of minimal
multiplicity ($2$ if $l=8$, $1$ otherwise). Then successive
contraction of $L,\d_l,...,\d_4,\d_3,\d_1$, in that order,
is the inverse of the blow-up described in (4), and
also proves (3).
\end{proof}
\end{proposition}

Now we can complete the proof of Theorem \ref{sst}.

At this point we know, by Proposition \ref{blank}
and Corollary \ref{configuration}, that $D_1^-$ is 
a weak del Pezzo surface of degree $d=9-l$ 
and that it contains a $(-2)$-configuration 
$\Delta=\{\delta_1,\ldots,\delta_l\}$ of type $G=E_l$.
So the configuration on $D_1^-$ formed by
the strict transform $\g^-$ of $\g$ and the flopped curves 
(see the paragraph following Lemma \ref{chain})
form part of the configuration $L,\delta_1,...,\delta_{l}$
given by Proposition \ref{6.30}
and so form a single $(-1,-2)$-chain. 
Therefore, by Lemma \ref{-1,-2 chain}, the flopping curves on $D_0$
form a unique $(-1,-2)$-chain on $D_0$.

This completes the proof of Theorem \ref{sst}.
\end{proof}
\end{theorem}
\end{section}
\begin{section}{Universal families of del Pezzo surfaces}\label{versal surfaces}
Theorem \ref{sst} refers to the line in the line bundle $\sL\to\sY$
that is the fibre over the origin $0_Y$, which is in turn 
defined over an algebraically closed
field. However, this extends over the whole of $\sL$ as follows.

Recall that $\sZ$ is obtained from a chart of $\sG_{\sE}$
by turning off $\Pic^1(\sE)$
and that $\tsZ=\sZ\times_{\sG_{\sE}}\tsG_{\sE}$.
There is a morphism $\rho:\tsZ\to\sL$ such that
the composite $\tsZ\to\sY$ is smooth
and $\rho^{-1}(0_{\sL})$ coincides with
the inverse image of $\Delta\times 0_{\sY}$,
where $\Delta$ is the discriminant divisor in $\sC_{pre}$, 
under the smooth
morphism $\tsZ\to\sC_{pre}\times_{\sM_{Ell}}\sY$.
So $\rho^{-1}(0_{\sL})=\sQ+\sD_0+\sD_1$.

Observe that there is a blowing-down morphism
$\tsZ\to\tsZ^+$
where $\sQ$ is contracted onto a curve in $\sD_1$,
because $\sQ$ is a family of quadrics with normal bundle
$\sO(-1,-1)$.

Recall the stack $\tsG^+$ from Section \ref{G+}.

\begin{lemma}\label{nnow} $\tsZ^+$ is isomorphic to 
$\sZ\times_{\sG}\tsG^+$.
\begin{proof} $\tsZ\to\tsZ^+$ is the contraction of a $\P^1\times\P^1$-bundle
$\sQ$ over $\sY$ to a $\P^1$-bundle $\g$. Fibre by fibre, this is the projection
$q:\P^1\times\P^1\to\P^1$ given by taking the family of $(1,1)$-curves
in $SL_3/B$ that pass through a fixed point and projecting to the family
of lines in $\P^2$ through a fixed point. Generically, these $(1,1)$-curves arise
as the rational tail $\phi$ in a pre-stable curve $C=E\cup\phi$.
From the construction of $\tsG^+$ and the morphism $\psi:\tsG\to\tsG^+$,
this projection is exactly achieved by $\psi$ over $\sZ$.
\end{proof}
\end{lemma}

\begin{proposition}\label{7.0}
\part[i] There are sequences of flops 
$\tsZ^+-\to\tsZ^-$
that extend the sequence of flops 
$\tX^+-\to\tX^-$.
\part[ii] They fit into a commutative diagram
$$\xymatrix{
&{\tsZ}\ar[dr]^{\rho}\ar[dl]\ar@{.>}[d]&&\\
{\tsZ^+}\ar@{.>}[r]\ar[dr]
&{\tsZ^-}\ar[r]_{\s}\ar[d]\ar@{}[dr]|-*+[F]{1}
&{\sL}\ar[d]\ar[r]&{\sY}\\
&{\sZ}\ar[r]&{[\hsY/W]}&
}$$
where the square marked
$\xymatrix@1{
*+[F]{1}
}$
is birationally Cartesian.
\part[iii] $\tsZ,\tsZ^+$ and $\tsZ^-$ are 
semi-stable over $\sL$. Their restrictions to $\sL-0_{\sL}$ are
isomorphic and smooth over $\sL-0_{\sL}$.
\part[iv] $\tsZ\vert_{0_{\sL}}=\sD_0+\sD_1+\sQ$,
$\tsZ^+\vert_{0_{\sL}}=\sD_0^++\sD_1^+$,
$\tsZ^-\vert_{0_{\sL}}=\sD_0^-+\sD_1^-$,
and $\sD_i$, $\sD_i^{\pm}$ and $\sQ$ are smooth over $\sY=0_{\sL}$.
\part[v] Each of $\sD_0\cap\sD_1$ and $\sD_0^{\pm}\cap\sD_1^{\pm}$ 
is isomorphic over $\sY$
to $\sE\times_{\sM_{Ell}}\sY$.
\part[vi] Each fibre of $\sD_0^-\to\sY$ is the minimal resolution
of the minimally elliptic singularity of degree $d$ that belongs
to the corresponding elliptic curve.
\begin{proof}
To prove the existence of the flops 
there is a local problem and a global problem to be solved.

The local problem is to show that under deformation 
of the ambient $3$-fold the
curves to be flopped do not disappear and that under
specialization they do not break up.

For each successive flop, the irreducible curve $\delta$ to
be flopped has normal bundle $\sN_{\delta}=\sO(-1)^{\oplus 2}$.
So $H^1(\delta,\sN_{\delta})=0$ and then $\delta$
survives under deformation.

After $\delta$ has been flopped, the number of irreducible curves
to be considered has diminished, and then induction shows that
all the curves to be flopped survive under deformation.

A similar argument shows that they do not break up under specialization.

The global problem
is to show that there is no monodromy 
acting on the lattice generated by the classes
$\b,\epsilon_1,...,\epsilon_t$ in the Chow group
of the generic fibre of the family $\sD_0\to\sY$.
This will follow from the next result.

\begin{proposition} $\sY$ is simply connected.
\begin{proof} We check first that $\sE$ is simply connected.

Suppose that $\pi:\tsE\to\sE$ is finite and {\'e}tale
and that $\tsE$ is connected.

Since $\sM_{Ell}$ is simply connected, by Theorem \ref{folklore},
the origin $0_{\sE}$ lifts to $\tsE$ and we have a $2$-commutative
diagram
$$\xymatrix{
{0_{\sE}}\ar@/^1pc/[drr]^i\ar[dr]^j\ar@/_1pc/[ddr]_{\cong}&&\\
&{\tsE}\ar[r]_{\pi}\ar[d]^{\alpha}&{\sE}\ar[dl]\\
&{\sM_{Ell}.}
}$$
Then $\b:\tsE\to\sM_{Ell}$ is an elliptic curve with 
identity $0_{\tsE}=j(0_{\sE})$.
This is induced by a morphism $\g:\sM_{Ell}\to\sM_{ell}$
and we get a commutative diagram
$$\xymatrix{
{0_{\sE}}\ar@/^1pc/[drr]^i\ar[dr]^j\ar@/_1pc/[ddr]_{\cong}&&\\
&{\sE\times_{\sM_{Ell},\a}\sM_{Ell}}\ar[r]_-{\pi}\ar[d]^{\b}&{\tsE}\ar[d]^{\alpha}\\
&{\sM_{Ell}}\ar[r]_{\g}&{\sM_{Ell}.}
}$$
Then $\g$ is an isomorphism, and so $\pi$ is an isomorphism
and $\sE$ is simply connected.

Since $\sY$ is a fibre product of copies of $\sE$
it is enough to check that, if $\sU\to\sM$ and $\sV\to\sM$
have sections $0_{\sU}$ and $0_{\sV}$ and if
$\sM,\sU$ and $\sV$ are simply connected, then so is
$\sU\times_{\sM}\sV=\sW$, say.

Suppose that $\tsW\to\sW$ is finite and {\'e}tale.
The sections $0_{\sU}$ and $0_{\sV}$ give sections
of $\sW\to\sU$ and of $\sW\to\sV$,
which then lift to sections $\sigma:\sU\to\tsW$
and $\tau:\sV\to\tsW$. This gives a section
$(\sigma,\tau):\sU\times_{\sM}\sV\to\tsW$,
as required.
\end{proof}
\end{proposition}
This completes the proof of Proposition \ref{7.0}.
\end{proof}
\end{proposition}

Corollary \ref{normal bundle} generalizes straightforwardly
to this situation.

\begin{proposition}
\part[i] Locally on $\sM_{Ell}$ there are $\sY$-isomorphisms $\tsZ^-\stackrel{\cong}{\to} 
N_{\sD^-_1/\tsZ^-}$ and $N_{\sD^-_1/\tsZ^-}\stackrel{\cong}{\to}\omega_{\sD^-_1/\sY}$.
\part[ii] Locally on $\sM_{Ell}$ there is an isomorphism,
relative to $\sY$, from
$\sZ$ 
to a family of $3$-fold cones, 
namely, to the cone over the anticanonical model of the family $\sD_1^-\to\sY$
of weak del Pezzo surfaces.
\noproof
\end{proposition}

Most of this section is devoted to proving that $\sD_1^-\to\sY$
is close to being universal.

Let $\tau:\tsZ^-\to\sD^-_1$ be the projection of
$N_{\sD^-_1/\tsZ^-}$ to its base. By construction,
there is a smooth morphism $\tsZ\to\tsG$;
this restricts to give a morphism
$$\tsZ^0:=\tsZ-\rho^{-1}(0_{\sL})=\tsZ^--\s^{-1}(0_{\sL})\to\tsG^{ss}.$$
We identify $\sD^-_1$ with the zero section of $N_{\sD^-_1/\tsZ^-}$
and $\sD^-_0$ with $\tau^{-1}(\sA)$. 

\begin{proposition}\label{descending lb}
There is an embedding $\Pic^G_F\inj\Pic_{\tsG^{ss}_{\sE}}$.
\begin{proof} A point $P$ of $\tsG^{ss}_{\sE}$ consists of
a $G$-bundle $\Xi$ over an elliptic curve $E$ and a section $\s$ of
$\sF_{\Xi}\to E$ whose cocharacter $[\s]$ vanishes. On the other hand,
an element $\varpi$ of $\Pic^G_F$ gives a line bundle $\sL_{\varpi}$ on
$\sF_{\Xi}$; evaluating $\sL_{\varpi}$ at the origin $0_\s$ of $\s$
gives a line, so a line bundle $\sM_{\varpi}$ on $\tsG^{ss}_{\sE}$
with $\sM_{\varpi}(P)=\sL_{\varpi}(0_\s)$.
This gives the embedding that was asserted. 
\end{proof}
\end{proposition}

Pull back to $\Pic(\tsZ^--\s^{-1}(0_{\sL}))$ via the morphism
$\tsZ^--\s^{-1}(0_{\sL})\to\tsG^{ss}$; we get
a homomorphism $\lambda:\Pic^G_F\to\Pic(\tsZ^--\s^{-1}(0_{\sL}))$.
Now fix an elliptic curve $E$ over a field $k$ and return to the restriction
$\tX^-\to\A^1_k$ of $\tsZ^-$ to the line $\A^1_k$ in $L$ over the origin
$0_Y$. For $t\ne 0$, the fibre $\tX^-_t$  
contains a $(-2)$-configuration $\D_t$
of type $G=E_l$; since $\tZ^-\cong N_{D^-_1/\tZ^-}$, $\D_t$ specializes
to an isomorphic such configuration $\D_0$ on $D^-_1$.

Regard $\tX^-_t$ as a surface contained in $\tsZ^--\s^{-1}(0_{\sL})$.
By Theorem \ref{classical}, $\D_t$ is identified with a subregular
unipotent Springer fibre associated to $G$ and, for any
$\varpi\in\Pic^G_F$, $\lambda(\varpi)$ is identified with $\varpi$.

We are led to the following result.

\begin{proposition} There are subgroups $\tH$ of $\Pic_{\tsZ}$ and
$H^-$ of $\Pic_{\tsZ^-}$ and a commutative diagram with exact rows
$$\xymatrix{
{\Z\{[\sD_0],[\sD_1],[\sQ]\}}\ar[r]\ar[d] & {\tH}\ar[r]\ar[d] & 
{\Pic^G_F}\ar[d]^{=}\ar[r] & {0}\\
{\Z\{[\sD_0^-],[\sD_1^-]\}}\ar[r] & {H^-} \ar[r] & {\Pic^G_F}\ar[r] & {0}
}$$
where the two leftmost vertical arrows are induced by the blowing-down
$\tsZ\to\tsZ^+$ and the flop $\tsZ^+-\to\tsZ^-$. 
In particular, $[\sQ]\mapsto 0$.
\begin{proof} The only remaining point is to check
that the obvious homomorphisms $\tH\to\Pic_{\tsZ}$ 
and $H^-\to \Pic_{\tsZ^-}$ are injective. However,
this follows from restricting to one of the surfaces
$\tZ_t$, as just described.
\end{proof}
\end{proposition}

Restrict from $\tsZ^-$ to $\sD_1^-$. We get a subgroup $H\subset\Pic_{\sD_1^-}$ and,
since $\sA=\sD_0^-\vert_{\sD_1^-}\sim \omega^{-1}_{\sD_1^-/\sY}$, an exact sequence
$$\xymatrix{
{\Z[\sA]\oplus\Z[\g]}\ar[r] & {H} \ar[r]^{\chi} & {\Pic^G_F}\ar[r] & {0}
}$$
such that, for any $\Lambda\in H$ and component $\d_i$ of $\D$, 
$(\Lambda.\d_i)_{\sD_1^-}=(\chi(\Lambda).\a_i^\vee),$
where $\a_i^\vee$ is the simple coroot corresponding to $\d_i$.

\begin{corollary} The natural homomorphism $H\to\NS(\sD_1^-/\sY)$ of sheaves
of commutative groups on $\sY$ is surjective, and $\NS(\sD_1^-/\sY)$ 
is constant.
\begin{proof} It is enough to prove surjectivity for one (geometric) weak
del Pezzo surface $D_1^-$. Choose $D_1^-$ lying over $0_Y$; this surface contains
a configuration $\D_0$ of type $G$ and a line $\g$, and the result now
follows from Proposition \ref{rigid dP} and the fact that $\sY$
is simply connected.
\end{proof}
\end{corollary}

Take $D_1^-$ as in the proof just given. Let $I_{1,l}$ denote the $\Z$-lattice
with $\Z$-basis $\{\d_1,...,\d_l,\g\}$ and inner product given by the intersection
numbers on $D_1^-$; then $I_{1,l}$ is isomorphic to the standard odd unimodular
hyperbolic lattice of rank $1+l$ and the basis just given describes an isometry
$\phi:I_{1,l}\to\NS(D_1^-)$. Under $\phi^{-1}$ the exceptional curves of the blow-up
$D_1\to\P^2$ described in Proposition \ref{rigid dP} are $\g,\g+\d_l,...,
\g+\d_l+\cdots+\d_3$ and $\g+\d_l+\cdots+\d_3+\d_1$, while $\d_2$ is the strict
transform of a line in $\P^2$. 

Since $\NS(\sD_1^-/\sY)$ is constant and $\Pic(\sD_1^-/\sY)\to \NS(\sD_1^-/\sY)$ is
an isomorphism, $\phi$ extends to an isometry $\phi:I_{1,l}\to\Pic(\sD_1^-/\sY)$.
Also, $\sD_1^-$ contains an anti-canonical divisor $\sA$ that is a copy of 
$\sE\times_{\sM_{Ell}}\sY$; the isomorphism 
$\sA\to \sE\times_{\sM_{Ell}}\sY$ is provided
by the base point $\g\cap\sA$.

Restricting to $\sA$ provides a homomorphism
$$\psi:I_{1,l}\to\Pic(\sA/\sY)\cong\Pic(\sE)=\coprod_{n\in\Z}\Pic^n(\sE).$$
Define $\kappa_l=\phi^{-1}[\omega^{-1}_{\sD_1/\sY}]$, the anti-canonical class.
Then 
$$\kappa_l=\d_1+\d_2+2\d_3+3\sum_{i\ge 4}\d_i+3\g.$$
Note that $I_{1,l}/\Z\kappa_l$ is the weight lattice $P$, while
$\oplus\Z\d_i$ is the root lattice $Q$. 
Consider the positive cone $C_+$ in $I_{1,l}\otimes\R$ defined by the conditions
$(v,v)\ge 0$, $v.\kappa_l\ge 0$;
this inherits a finite decomposition
into chambers from the decomposition of $P\otimes\R$ into fundamental domains
for the $W$-action. One such is the chamber $C_0$
defined by the inequalities $(v,\d_i)\ge 0$ for all $i$.

\begin{lemma}\label{degree}
$\omega^{-1}_{\sD_1^-/\sY}\vert_{\sA}$
is isomorphic to the pullback of $\sO_E(d[O_E])$
to $\sA=\sE\times_{\sM_{Ell}}\sY$.
\begin{proof}
$\sD_0^-$ is the minimal resolution of a singular normal surface
that is defined over the stack $\sM_{Ell}$ of elliptic curves,
and so $\sN_{\sA/\sD_0^-}$ is defined over $\sM_{Ell}$.
So $\sN_{\sA/\sD_0^-}\cong\sO_E(-d[O_E])$.
Since $\sN_{\sA/\sD_0^-}\cong \sN_{\sA/\sD_1^-}^\vee$,
the result follows from the adjunction formula
and the triviality of $\omega_{\sZ^-/\sL}$.
\end{proof}
\end{lemma}

By Lemma \ref{degree},
$\psi(\kappa_l)=d[O_E]$,
so there is a homomorphism
$\barpsi: P\to \Pic^0_E$ 
defined by $\barpsi(L)=\psi(L)-\deg\psi(L)[O_E]$.

Now consider the stack $\sM d\sP_d$ of {\emph{marked}} weak
del Pezzo surfaces of degree $d$; by definition, the objects consist of:

\noindent (1) a family $f:X\to S$ of weak del Pezzo surfaces of degree $d$;

\noindent (2) an embedding $i:E\times S\inj X$ such that the 
class of the image
$\sA$ equals the class $[\omega^{-1}_{X/S}]$ and $i^*(\omega^{-1}_{X/S})$
is linearly equivalent to $d[O_E]$;

\noindent (3) an isometry $\phi:I_{1,l}\to \Pic(X/S)$ such that
$\phi(\g)\vert_{\sA}$ is linearly equivalent to $[0_{\sA}]$ and
$\phi(\d_j)$ has degree $0$ for all $j$.

In this language, $(\sD_1^-\to\sY,\ \sA\inj\sD_1^-,\ \phi)$
is an object of $\sM d\sP_d$, so defines a morphism $H:\sY\to\sM d\sP_d$.

There is also a morphism $F:\sM d\sP\to\sY$ defined by
$(f,i,\phi)\mapsto\psi'$, where $\psi'$ is constructed exactly as above, 
and a morphism $G:\sY\to \sM d\sP$, as follows:

Given $\psi':P\to E$, construct $\psi:I_{1,l}\to\Pic(E)$ by 
$$\psi(\g)=[0_{E}],\ \psi(\d_i)=\psi'(\d_i{\pmod{\Z\kappa_l}}).$$
Then $\psi(\kappa_l)$ has degree $3$. Embed $E\inj\P^2$
via $\vert\psi(\kappa_l)\vert$ and then make $l$ blow-ups $\P^2$
along the points on $E$ defined by the degree $1$ classes
$$\psi(\g+\d_l+\cdots+\d_3+\d_1),\psi(\g+\d_l+\cdots+\d_3),...,\psi(\g+\d_l),\psi(\g)$$
in that order.

The next result is due to M{\'e}rindol [M]. Its proof is
an immediate consequence of the definitions of $F$ and $G$.

\begin{proposition}
$F\circ G=1_{\sY}$.
\noproof
\end{proposition}

\begin{corollary} \label{sep quotient} $F$ exhibits $\sY$ as the 
maximal separated quotient of $\sM d\sP$.
\begin{proof}
Immediate.
\end{proof}
\end{corollary}

That is, $\sM d\sP$ is obtained by glueing together copies of $\sY$ along
open subvarieties. On the locus of del Pezzo surfaces that have no $(-2)$-curves
the stacks $\sM d\sP$ and $\sY$ are isomorphic but $\sM d\sP$ fails
to be separated when $(-2)$ curves appear.

\begin{proposition} $F\circ H=1_{\sY}$.
\begin{proof} $F\circ H$ takes $0_{\sY}$ to $0_{\sY}$, so is a homomorphism
of abelian schemes. Fix a geometric elliptic curve $E$; then over $0_Y$
the fibre of $\sD_1$ contains a $(-2)$-configuration of type $G$. 
Consider the subdiagram
$$
\xymatrix{
{\tsZ^-}\ar[r]^{\s}\ar[d]_{\epsilon} & {\sL}\ar[r]\ar[d] & {\sY}\\
{\sZ}\ar[r] & {[\hsY/W].}
}$$
of the diagram in Proposition \ref{7.0}.
It follows from Theorem \ref{HS} that $\sZ\to [\hsY/W]$, which is a family 
of affine surfaces, has du Val singularities of type $G$ exactly over
$\A^1-\{0\}$, where $\A^1$ is the line in the cone $[\hsY/W]$ that lies over
the image of $0_{\sY}$ in $[\hsY/W]$. Therefore $0_{\sY}$ is isolated
in the fibre $(F\circ H)^{-1}(0_{\sY})$ and then $F\circ H$ is an isomorphism.
Since $F\circ H$ is $W$-equivariant it is then $\pm 1$.

Suppose that $F\circ H=-1$. Put $S=\sU\times_{\sM d\sP, G}\sY$, where
$\sU\to\sM d\sP$ is universal, and consider the $\sY$-isomorphism
$\sD_1\to S\times_{\sY,F\circ H}\sY$. Since $F\circ H=-1$, this takes
the ample cone on $\sD_1$ to the negative of the ample cone on $S$.
This is impossible, so $F\circ H=1$.
\end{proof}
\end{proposition}

The modular interpretation of $\sY$ as an open substack of $\sM d\sP$
is given as follows.

Suppose that $X$ is a weak del Pezzo surface. Define a \emph{line}
in $X$ to be a class $m$ in $\NS(X)$ such that 
$m^2=-1$ and $m.(-K_X)=1$ and a \emph{root} in $X$ to be
a class $\d$ in $\NS(X)$ such that 
$\d^2=-2$ and $\d.K_X=0$. Let $\Lambda$ denote the set of lines
and $\Delta$ the set of roots on $X$. Then $\Delta$ is a finite root system
and reflections in the roots generate a finite Weyl group
$W$, as usual. A root $\d$ is \emph{effective}
if $H^0(X,\d)\ne 0$.

Define the \emph{positive cone}
$C^+=C^+(X)$ by $C^+=\{x\in\NS(X)_\R\mid x^2>0,\ x.(-K_X)>0\}$.
It is well known that a class $x$ in $C^+\cap\NS(X)$ is nef if and only
if $x.m>0$ for all $m\in\Lambda$ and $x.\d\ge 0$ for all effective
roots $\d$ on $X$.
Put $\sF=\{x\in C^+\mid x.m>0\ \forall m\in\Lambda\}$.
Then the roots define walls that tesselate $\sF$ into chambers
that are permuted simply transitively by $W$.
So, if $\sF'$ is one such chamber, then $\sF=\cup_{w\in W}w(\sF')$.

In the lattice $I_{1,l}$, say 
$\Lambda_0=\{m\in I_{1,l}\mid m^2=-1,\ m.\kappa =1\}$,
$\Delta_0=\{\d\in I_{1,l}\mid \d^2=-2,\ \d.\kappa=0\}$,
$C_0^+=\{x\in I_{1,l}\otimes\R\mid x^2>0,\ x.\kappa>0\}$
and $\sF_0=\{x\in C_0^+\mid x.m>0\ \forall m\in\Lambda_0\}$.

Pick a chamber $\sF_0'$ in the tessellation of 
$\sF_0$ defined by the roots in $\Delta_0$.
Let $\sM d\sP^+$ be the 
open substack of $\sM d\sP$ whose objects
are triples $(X\to S,i,\phi)$ such that
$\phi_\R(\sF_0')$ is contained in the nef cone of every
geometric fibre of $X\to S$. Then $\sM d\sP^+$ 
is isomorphic to $\sY$, and $\sM d\sP$
is the union of these copies of $\sY$, one for each
chamber $\sF_0'$ in $\sF_0$.

This brings us to the main result of the paper.

\begin{theorem}\label{main} 
Suppose that $G=E_l$ for $l=5,6,7,8$.

\part[i] The family 
$\sZ\to[\hsY/W]$ has simply elliptic singularities
over the vertex of the cone $[\hsY/W]$.

\part[ii] Over $\Sp\Z[1/(9-l)]$ this family
is a miniversal deformation
of these singularities. 

\part[iii] The commutative diagram
$$\xymatrix{
{\tsZ^-}\ar[r]\ar[d] & {\sZ}\ar[d]\\
{\sL}\ar[r] & {[\hsY/W]}
}$$
is a simultaneous log resolution of $\sZ\to [\hsY/W]$.

\part[iv] The exceptional divisor $\sD_1^-$ in
$\tsZ^-$ is the 
restriction to $\sY$ of the universal marked del Pezzo
surface of degree $d=9-l$.
\begin{proof}
\DHrefpart{i}-\DHrefpart{iii} summarize the results
of Section \ref{subregular}.
\DHrefpart{iv}
follows at once from Corollary \ref{sep quotient}
and the construction of $\tsZ^-$.
\end{proof}
\end{theorem}

\begin{remark}
Suppose that the base is $\Sp\C$.
The affine del Pezzo surface obtained by deleting the elliptic curve $E$
has a mixed Hodge structure on $H^2$. These mixed Hodge structures are 
naturally parametrized by the Looijenga variety $\sY$. The above theorem
shows, when restricted to the exceptional divisor $\sD_1^-$ in $\sZ^-$,
that the period map for these affine del Pezzo surfaces is the same as
the map from $\sD_1^-$ to $\sY$ that we have constructed in terms of 
group theory.
\end{remark}
\end{section}
\begin{section}{del Pezzo surfaces
and the unipotent singularity of $E_8$ 
in characteristics $2,3$ and $5$}\label{weird}
Consider again the family $\tX^-\to\A^1$
where $\A^1$ is the line in the line bundle $L\to Y$
lying over $0_Y$.
The closed fibre of $\tX^-_0$ is $\tX^-_0=D_0^-+D_1^-$
and for $t\ne 0$ the fibre $\tX^-_t=\tX_t$
contains a configuration $\G_t$
of $(-2)$-curves of the same combinatorial type as $G$.
When $t=0$ this configuration specializes to an
isomorphic configuration
$\G_0$ which lies in $D_1^-$ and
is disjoint from the double curve
$A^-=D_0^-\cap D_1^-$. Recall that $A^-$
is isomorphic to $E$.

Let $D_1^-\to D_1^\flat$ be the contraction of $\G_0$.
Then $D_1^\flat$ has a singularity of combinatorial type $G$
that does not lie on the curve $A^\flat$ which is the
(isomorphic) image of $A^-$.

\begin{proposition}\label{preceding} 
\part[i] $\tX^-$ is isomorphic to the 
line bundle $\omega_{D_1^-}$.
\part[ii] There is a commutative diagram
$$\xymatrix{
{\tX^-}\ar[r]\ar[dr]&{\tX^\flat}\ar[d]\\
&{\A^1_k}
}$$
where $\tX^-\to \tX^\flat$ contracts each $\G_t$ to a du Val singularity,
$\cup_t\G_t$ is contracted to a section
of $\tX^\flat\to\A^1_k$ and $\tX^\flat_0=D_0^\flat+D_1^\flat$
where $D_0^+\to D_0^\flat$ is an isomorphism.
\part[iii] $\tX^\flat$ is isomorphic to the line bundle
$\omega_{D_1^\flat}$.
\part[iv]
The singular affine surface $\tX_t$
is isomorphic to the affine del Pezzo surface
$D_1^\flat-A^\flat$ for all $t\ne 0$.
\begin{proof} The existence of the contraction $X^+\to X^\flat$
is a special case of the well known
fact that a family of $(-2)$-configurations
lying in the relative smooth locus of
a family of surfaces can be simultaneously
contracted.
The rest follows from Lemma \ref{6.10}.
\end{proof}
\end{proposition}

As has been aready recalled,
if $Uni(G)$ is the unipotent variety of a 
simply connected simple
group $G$ of type $A,D$ or $E$, then the singular
locus of $Uni(G)$ is the subregular locus of $Uni(G)$
and the strict localization of $Uni(G)$ 
at the geometric generic point of its subregular locus
(we shall refer to this as the unipotent singularity)
has a du Val singularity
of the same combinatorial type as $G$. However, in low characteristics
the combinatorial type does not specify the singularity,
even over an algebraically closed field;
for example, the $E_8$ singularities are classified up to 
formal isomorphism 
on p. 270 of [Li] (and up to henselian equivalence in [Ar]).
The local equation of each is
$$t^2+z^3+y^5=\phi$$
where $\phi\in\{0,y^3t,y^3zt,y^2zt,yzt\}$ in characteristic $2$,
$\phi\in \{0,y^3z^2,y^2z^2\}$ in characteristic $3$ and
$\phi\in\{0,y^4z\}$ in characteristic $5$.
We refer to them as $E_8^{(p,\phi)}$ accordingly.
In Artin's list [Ar] they are listed instead as
$E_8^r$ where $0\le r\le 4$ or $2$ or $1$ in these characteristics
and
$\dim T^1=16-2r$ or $12-2r$ or $10-2r$.
 
For the rest of this section
the base field $k$ will be
a complete valued algebraically closed
field of characteristic $p>0$ and $E$ will be a uniformizable
elliptic curve over $k$.

\begin{theorem} 
The singularity on $D_1^\flat$
is isomorphic to the unipotent singularity of
the group $G$.
\begin{proof} By Theorem \ref{Springer}
the singularity of $X_{t\ne 0}$
is isomorphic to the unipotent singularity.
The theorem follows from Proposition \ref{preceding}.
\end{proof}
\end{theorem}

Recall that in characteristic $p$ the $j$-invariant of
a uniformizable elliptic curve $E$ is not algebraic over
$\F_p$; in particular, it is non-zero.

Until further notice the group $G$ will be $E_8$.

Recall that the anticanonical model of a del Pezzo surface $S$ of degree $1$
over a field $k$ 
is a sextic hypersurface in weighted projective space 
$\P(1,1,2,3)=\Proj k[X_1,Y_1,Z_2,T_3]$. 

The next result is a counterpart to Proposition \ref{rigid dP}.

\begin{proposition}\label{del Pezzo}
Suppose that $S$ is a del Pezzo surface of degree $1$ 
over a separably closed field $k$ with
a singular point $P$ of type $E_8$.

\part [i] Up to projective equivalence,
the defining equation of $S$ in $\P(1,1,2,3)$ is
$$T^2+Z^3+XY^5=\Phi$$ where
$\Phi\in\{ZY^4,0\}$ if $\ch k\ne 2,3$,
$\Phi\in\{Z^2Y^2,ZY^4,0\}$ if $\ch k=3$ and
$\Phi\in\{TYZ,TY^3,\g ZY^4\}$ if $\ch k=2$,
where $\g=\lambda^4+c$, so defines a class in $H^1_{flat}(k,\mu_4)$.
\part[ii] Suppose that $C$
is a smooth member of $\vert-K_S\vert$. 
Then, whatever the characteristic, $j(C)$ is never zero for
the first value of $\Phi$
and is always zero for the second. (In characteristics $2$ and $3$
the pencil $\vert-K_S\vert$ is quasi--elliptic for the third value, namely,
when $\Phi=0$ or $\g ZY^4$.)
\part[iii] 
In the first case of each characteristic the singularity $(S,P)$ is isomorphic to 
$E_8^{(5,y^4z)}=E_8^1$ or $E_8^{(3,y^2z^2)}=E_8^2$ or $E_8^{(2,yzt)}=E_8^4$.
\begin{proof} This depends upon the following recognition principle
for $E_8$ singularities. It goes back to du Val.
\begin{lemma}
\part[i] A du Val singularity $(S,P)$ is 
of combinatorial type $E_8$
if and only if 
the blow-up $\tS=\Bl_PS$ has a singular point
of type $E$.
\part[ii] A du Val singularity 
whose defining equation can be written as
$$t^2-G(y,z)\in t(y,z,t)^2$$
is of type $E$ if and only if the leading term of $G$
is proportional to a perfect cube.
\begin{proof}
\DHrefpart{i} Suppose that $S^*\to S$ is the minimal resolution.
Then, according to du Val and Artin,
this factors through $\tS$ and the fundamental
cycle $Z$ on $S^*$ is a divisor and the sheaf of ideals
that defines it is
$\frak m_P.\sO_{S^*}$. So, for any component $C_i$
of the exceptional locus in $S^*$, $Z.C_i\le 0$,
and $Z.C_i<0$ if and only if $C_i$ appears on $\tS$.
Moreover, $Z$ is the biggest
root; in the case of $E_8$ this is
the following linear combination of simple roots:
$$\xymatrix{&&{3}\ar@{-}[d]&&&&\\       
{2}\ar@{-}[r]&{4}\ar@{-}[r]&{6}\ar@{-}[r]&{5}\ar@{-}[r]       
&{4}\ar@{-}[r]&{3}\ar@{-}[r]&{2}                                      
}$$
Therefore the curves contracted by $S^*\to\tS$
form an $E_7$ configuration.
Similarly an $E_7$ (or $E_6$) singularity has a $D_6$ 
(or $A_5$) singularity on its
first blow-up.

\DHrefpart{ii} This is to be found on p. 268 of [Li]. 
\end{proof}
\end{lemma}
Beyond this lemma the proof of Proposition \ref{del Pezzo} consists of a
straightforward manipulation of polynomials
and a calculation of $j$-invariants using the
standard formulae to be found on, for example,
p. $46$ of [Si].
We omit the details.

\end{proof}
\end{proposition}
\begin{remark} In particular, over an algebraically
closed field of characteristic $2$ there are
five $E_8$ singularities but only three of them can lie on
a del Pezzo surface. We do not know what to make of this.
\end{remark}
\begin{corollary}\label{unipotent sing}
The unipotent singularity of the group $E_8$
in characteristics $2,3$ and $5$ is isomorphic to
$E_8^{(2,yzt)}=E_8^4$, $E_8^{(3,y^2z^2)}=E_8^2$ and
$E_8^{(5,y^4z)}=E_8^1$, respectively.
\begin{proof}
The unipotent singularity is realized on a del Pezzo
surface that contains a uniformizable curve.
The $j$-invariant of a uniformizable curve $E$ is 
transcendental over $\F_p$,
so that in particular $j(E)\ne 0$, and we can
now apply Proposition \ref{del Pezzo}.
\end{proof}
\end{corollary}
\begin{corollary} The BGSS construction gives
a versal deformation of the unipotent singularity of $E_8$
in all characteristics.
\begin{proof}
For each du Val singularity Artin computed [Ar] the dimension of $T^1$,
the tangent space to a miniversal deformation space
of the singularity.
From his list and Corollary \ref{unipotent sing} we can see that
$\dim T^1=8$ for the unipotent singularity of $E_8$.
Now the argument of [SB01] goes through unchanged,
as it does whenever $\dim T^1$ equals the rank of the group.
\end{proof}
\end{corollary} 
\begin{corollary}\label{Springer false} 
If $E$ is supersingular, $G=E_8$ and $\ch k=2,3$ or $5$
then $\tsG^{ss}\to\sG^{ss}$ is not smoothly equivalent
to $\tG\to G$ in any neighbourhood of the identity of $G$.
That is, the conclusion of Theorem \ref{Springer}
fails in these cases.
\begin{proof} If it were true, then
the morphism $\tsG^{ss}\to\sG^{ss}$ would
be smoothly equivalent to $\tG\to G$
and then the unipotent singularity of $E_8$ would be realized
on a del Pezzo surface $S$
on which $E$ is an anticanonical curve.
Then, according to
Proposition \ref{del Pezzo}, $S$ can be defined by the equation
$$T^2+Z^3+XY^5=\Phi$$
where $\Phi=0,ZY^4$ or $TY^3$ when $\ch k=5,3$ or $2$,
respectively,
since in these characteristics an elliptic curve is supersingular
if and only if its $j$-invariant is $0$.
But then $S$ has a singularity of type
$E_8^{(5,0)}$ or $E_8^{(3,y^4z)}$ or $E_8^{(2,y^3t)}$
and it remains only to check that $E_8^{(3,y^4z)}$ is
not isomorphic to $E_8^{(3,y^2z^2)}$. This can be done by
calculating $\dim T^1$ for $E_8^{(3,y^4z)}$;
the result is $\dim T^1=10$, so that inspection of Artin's list [Ar]
shows that $E_8^{(3,y^4z)}$ is in fact isomorphic to $E_8^{(3,y^3z^2)}$.
\end{proof}
\end{corollary}
For $E_6$ and $E_7$ similar considerations involving
del Pezzo surfaces of degrees $3$ and $2$ work to
describe the unipotent singularity in low characteristics:
in the notation of [Ar]
they are $E_6^1$ and $E_7^3$ when $p=2$ and
$E_6^1$ and $E_7^1$ when $p=3$.
In all cases the unipotent singularity is the one in 
its combinatorial class where $\dim T^1$ is minimal,
so is the ``most general'' one. 
\end{section}
\bibliography{alggeom,ekedahl}
\bibliographystyle{pretex}
\end{document}